%% file: extremII.tex
\documentstyle[11pt,amsfonts,fullpage]{amsart}

\def\endproof{\qed\medskip}
\def\blacksquare{\hbox to .60em{\vrule width .60em height .60em}}

\newcommand{\Roof}[2]{#2{#1}}

\newcommand{\bbgin}{\begin}

\begin{document}

\title[ ]{Extrema of Curvature Functionals on the Space of Metrics on 3-Manifolds, II.}

\author[ ]{Michael T. Anderson}

\thanks{Partially supported by NSF Grant DMS-9802722}

\maketitle

\setcounter{section}{-1}

\section{Introduction}
\setcounter{equation}{0}

 This paper is a continuation of the study of some rigidity or non-existence issues discussed in [An1, \S 6]. The results obtained here also play a significant role in the approach to geometrization of 3-manifolds discussed in [An4].

 Let $N$ be an oriented 3-manifold and consider the functional
\begin{equation} \label{e0.1}
{\cal R}^{2}(g) = \int_{N}|r_{g}|^{2}dV_{g}, 
\end{equation}
on the space of metrics ${\Bbb M} $ on $N$ where $r$ is the Ricci curvature and $dV$ is the volume form. The Euler-Lagrange equations for a critical point of ${\cal R}^{2}$ read
\begin{equation} \label{e0.2}
\nabla{\cal R}^{2} = D^{*}Dr + D^{2}s -  2 \stackrel{\circ}{R}\circ r -\tfrac{1}{2}(\Delta s -  |r|^{2})\cdot  g = 0, 
\end{equation}

\begin{equation} \label{e0.3}
\Delta s = -\tfrac{1}{3}|r|^{2}. 
\end{equation}
Here $s$ is the scalar curvature, $D^{2}s$ the Hessian of $s$, $\Delta s = trD^{2}s$ the Laplacian, and $ \stackrel{\circ}{R} $ the action of the curvature tensor $R$ on symmetric bilinear forms, c.f. [B, Ch.4H] for further details. The equation (0.3) is just the trace of (0.2). 

 It is obvious from the trace equation (0.3) that there are no non-flat ${\cal R}^{2}$ critical metrics, i.e. solutions of (0.2)-(0.3), on compact manifolds $N$; this follows immediately by integrating (0.3) over $N$. Equivalently, since the functional ${\cal R}^{2}$ is not scale invariant in dimension 3, there are no critical metrics $g$ with ${\cal R}^{2}(g) \neq $ 0. To obtain non-trivial critical metrics in this case, one needs to modify ${\cal R}^{2}$ so that it is scale-invariant, i.e. consider $v^{1/3}{\cal R}^{2},$ where $v$ is the volume of $(N, g)$.

 Nevertheless, it is of course apriori possible that there are non-trivial solutions of (0.2)-(0.3) on non-compact manifolds $N$.
\bbgin{theorem} \label{t 0.1.}
  Let (N, g) be a complete ${\cal R}^{2}$ critical metric with non-negative scalar curvature. Then (N,g) is flat.
\end{theorem}

 This result generalizes [An1, Thm.6.2], which required that $(N, g)$ have an isometric free $S^{1}$ action. It is not known if the condition $s \geq $ 0 is necessary in Theorem 0.1; a partial result without this assumption is given after Proposition 2.2. However, following the discussion in \S 1 and [An1, \S 6], the main situation of interest is when $s \geq $ 0.

 Of course, Theorem 0.1 is false in higher dimensions, since any Ricci-flat metric is a critical point, in fact minimizer of ${\cal R}^{2}$ in any dimension, while any Einstein metric is critical for ${\cal R}^{2}$ in dimension 4. 

 Next, we consider a class of metrics which are critical points of the functional ${\cal R}^{2}$ subject to a scalar curvature constraint. More precisely, consider scalar-flat metrics on a (non-compact) 3-manifold $N$ satisfying the equations
\begin{equation} \label{e0.4}
\alpha\nabla{\cal R}^{2} + L^{*}(\omega ) = 0, 
\end{equation}

\begin{equation} \label{e0.5}
\Delta\omega  = -\frac{\alpha}{4}|r|^{2}. 
\end{equation}
where again (0.5) is the trace of (0.4) since $s =$ 0. Here $L^{*}$ is the adjoint of the linearization of the scalar curvature, given by
$$L^{*}f = D^{2}f -  \Delta f\cdot  g - fr, $$
and $\omega $ is a locally bounded function on $N$, which we consider as a potential. The meaning and derivation of these equations will be discussed in more detail in \S 1. They basically arise from the Euler-Lagrange equations for a critical metric of ${\cal R}^{2}$ subject to the constraint $s =$ 0.

 The parameter $\alpha $ may assume any value in [0, $\infty ).$ When $\alpha  =$ 0, the equations (0.4)-(0.5) are the static vacuum Einstein equations, c.f. [An2] and references there. In this case, we require that $\omega $  is not identically 0. 

 It is proved in \S 3 that an $L^{2,2}$ Riemannian metric $g$ and $L^{2}$ potential function $\omega $ satisfying the equations (0.4)-(0.5) weakly in a 3-dimensional domain is a $C^{\infty},$ (in fact real-analytic), solution of the equations. A smooth metric $g$ and potential function $\omega $ satisfying (0.4)-(0.5) will be called an ${\cal R}_{s}^{2}$ critical metric or ${\cal R}_{s}^{2}$ solution.
\bbgin{theorem} \label{t 0.2.}
  Let $(N, g)$ be a complete ${\cal R}_{s}^{2}$ critical metric, i.e. a complete scalar-flat metric satisfying (0.4)-(0.5), with 
\begin{equation} \label{e0.6}
-\lambda  \leq  \omega  \leq  0, 
\end{equation}
for some $\lambda  <  \infty .$ Suppose further that (N, g) admits an isometric free $S^{1}$ action leaving $\omega $ invariant. Then (N, g) is flat.
\end{theorem}

 In contrast to Theorem 0.1, the assumption that $(N, g)$ admit an isometric free $S^{1}$ action here is essential. There are complete non-flat ${\cal R}_{s}^{2}$ solutions satisfying (0.6) which admit an isometric, but not free, $S^{1}$ action. For example, the complete Schwarzschild metric is an ${\cal R}_{s}^{2}$ solution for a suitable choice of the potential $\omega $ satisfying (0.6), c.f. Proposition 5.1.

 The condition $\omega  \leq $ 0 can be weakened to an assumption that $\omega  \leq $ 0 outside some compact set $K \subset N$. However, it is unknown if this result holds when $\omega  \geq $ 0 everywhere for instance. Similarly, the assumption that $\omega $ is bounded below can be removed in certain situations, but it is not clear if it can be removed in general, c.f. Remarks 4.4 and 4.5.

 The proofs of these results rely almost exclusively on the respective trace equations (0.3) and (0.5). The full equations (0.2) and (0.4) are used only to obtain regularity estimates of the metric in terms of the potential function $s$, respectively $\omega .$ Thus it is likely that these results can be generalized to variational problems for other curvature-type integrals, whose trace equations have a similar form; c.f. \S 5.2 for an example.

 As noted above, both Theorem 0.1 and 0.2 play an important role in the approach to geometrization of 3-manifolds studied in [An4]. For instance, Theorem 0.2 is important in understanding the collapse situation.

 Following discussion of the origin of the ${\cal R}^{2}$ and ${\cal R}_{s}^{2}$ equations in \S 1, Theorem 0.1 is proved in \S 2. In \S 3, we prove the regularity of ${\cal R}_{s}^{2}$ solutions and apriori estimates for families of such solutions. Theorem 0.2 is proved in \S 4, while \S 5 shows that the Schwarzschild metric is an ${\cal R}_{s}^{2}$ solution and concludes by showing that Theorems 0.1 and 0.2 also hold for ${\cal Z}^{2}$ and ${\cal Z}_{s}^{2}$ solutions, where $z$ is the trace-free Ricci curvature. 

  While efforts have been made to make the paper self-contained, in certain instances we refer to [An1] for further details.

\section{Scalar Curvature Constrained Equations.}
\setcounter{equation}{0}

 In this section, we discuss the nature and form of the ${\cal R}_{s}^{2}$ equations, as well as some motivation for considering these and the ${\cal R}^2$ equations. The discussion here is by and large only formal and we refer to [An1, \S 8] and [An4] for complete details and proofs of the assertions made.

 Suppose first that $M$ is a compact, oriented 3-manifold and consider the scale-invariant functional
\begin{equation} \label{e1.1}
I_{\varepsilon} = \varepsilon v^{1/3}{\cal R}^{2} + v^{1/6}{\cal S}^{2} = \varepsilon v^{1/3}\int|r|^{2} + v^{1/6}(\int s^{2})^{1/2}. 
\end{equation}
on the space of metrics on $M$. Here $\varepsilon  > $ 0 is a small parameter and we are interested in considering the behavior $\varepsilon  \rightarrow $ 0. The existence, regularity and general geometric properties of minimizers $g_{\varepsilon}$ of (essentially) $I_{\varepsilon}$, for a fixed $\varepsilon > 0$, are proved in [An1, \S 8]; more precisely, such is done there for the closely related functional $\varepsilon v^{1/3}\int|R|^{2} + v^{1/3}\int s^{2}$, where $R$ is the full Riemann curvature tensor. All of these results follow from the same results proved in [An1, \S 3 - \S 5] for the $L^2$ norm of $R$, together with the fact that, for a fixed $\varepsilon > 0$, the functional $\varepsilon v^{1/3}\int|R|^{2} + v^{1/3}\int s^{2}$ has the same basic properties as $v^{1/3}\int|R|^{2}$ w.r.t. existence, regularity and completeness issues. Now in dimension 3, the full curvature $R$ is controlled by the Ricci curvature $r$. Thus, for example one has the relations $|R|^2 = 4|r|^2 - s^2$ and $s^2 \leq 3|r|^2$. A brief inspection of the work in [An1, \S 3 - \S 5, \S 8] then shows that these results, together with the same proofs, also hold for the functional $I_{\varepsilon}$ and its minimizers $g_{\varepsilon}$ in (1.1).

  The Euler-Lagrange equations for $I_{\varepsilon}$ at $g = g_{\varepsilon}$ are
\begin{equation} \label{e1.2}
\varepsilon\nabla{\cal R}^{2} + L^{*}(\tau ) + (\tfrac{1}{4}s\tau  + c)\cdot  g = 0, 
\end{equation}
\begin{equation} \label{e1.3}
2\Delta (\tau + {\tfrac{3}{4}}\varepsilon s) + {\tfrac{1}{4}}s\tau  = -{\tfrac{1}{2}}\varepsilon |r|^{2} + 3c, 
\end{equation}
where $\tau  = \tau_{\varepsilon} \ s/\sigma , \sigma  = (v^{1/6}{\cal S}^{2}(g))$ and the constant term $c$, corresponding to the volume terms in (1.1), is given by
\begin{equation} \label{e1.4}
c = \frac{\varepsilon}{6v}\int|r|^{2}dV + \frac{1}{12\sigma v}\int s^{2}dV. 
\end{equation}
Again (1.3) is the trace of (1.2). These equations can be deduced either from [An1, \S 8] or [B, Ch.4H]; again all terms in (1.2)-(1.4) are w.r.t. $g = g_{\varepsilon}$. 

 As $\varepsilon  \rightarrow $ 0, the curvature $r_{\varepsilon}$ of the solutions $g_{\varepsilon}$ of (1.2)-(1.3) will usually blow-up, i.e. diverge to infinity in some region, say in a neighborhood of points $x_{\varepsilon}.$ Thus blow up or renormalize the metric $g_{\varepsilon}$ by considering $g_{\varepsilon}'  = \rho_{\varepsilon}^{- 2}\cdot  g_{\varepsilon},$ where $\rho_{\varepsilon} \rightarrow $ 0 is chosen so that the curvature in the geodesic ball $(B_{x_{\varepsilon}}' (1), g_{\varepsilon}' )$ w.r.t. $g_{\varepsilon}' $ is bounded. More precisely, $\rho $ is chosen to be $L^{2}$ curvature radius of $g_{\varepsilon}$ at $x_{\varepsilon},$ c.f. [An1,Def.3.2]. Then the renormalized Euler-Lagrange equations take the form
\begin{equation} \label{e1.5}
\frac{\varepsilon}{\rho^{2}}\nabla{\cal R}^{2} + L^{*}(\tau ) +({\tfrac{1}{4}}s\tau  + \frac{c}{\rho^{2}})\cdot  g = 0,  
\end{equation}

\begin{equation} \label{e1.6}
2\Delta (\tau +{\tfrac{3}{4}}\frac{\varepsilon}{\rho^{2}}s) + {\tfrac{1}{4}}s\tau  = - {\tfrac{1}{2}}\frac{\varepsilon}{\rho^{2}}|r|^{2} + \frac{3c}{\rho^{2}}. 
\end{equation}
Here $\rho  = \rho_{\varepsilon}$ and otherwise all metric quantities are taken w.r.t. $g = g_{\varepsilon}' $ except for the potential function $\tau ,$ which has been normalized to be scale invariant. It is necessary to divide by $\rho^{2}$ in (1.5)-(1.6), since otherwise all terms in the equations tend uniformly to 0 in $L^{2}(B_{\varepsilon}), B_{\varepsilon} = B_{x_{\varepsilon}}' (1).$

 From the scaling properties of $c$ in (1.4), note that in the scale $g_{\varepsilon}' , c'  = \rho^{4}\cdot  c,$ where $c = c(g_{\varepsilon}).$ Thus, the constant term in (1.5)-(1.6) satisfies
\begin{equation} \label{e1.7}
\frac{c}{\rho^{2}} = \frac{c(g_{\varepsilon}' )}{\rho^{2}} = \frac{\rho^{4}c(g_{\varepsilon})}{\rho^{2}} = \rho^{2}c(g_{\varepsilon}) \rightarrow  0, \ \ {\rm as} \ \ \varepsilon  \rightarrow  0. 
\end{equation}
Similarly, since ${\cal S}^{2}(g_{\varepsilon})$ is bounded, and the scalar curvature $s = s_{\varepsilon}' $ of $g_{\varepsilon}' $ is given by $s_{\varepsilon}'  = \rho^{2}s_{\varepsilon},$ one sees that $s$ in (1.5)-(1.6) goes to 0 in $L^{2}(B_{\varepsilon})$ as $\varepsilon  \rightarrow $ 0.

 Now assume that the potential function $\tau = \tau_{\varepsilon}$ is uniformly bounded in the $g_{\varepsilon}' $ ball $B_{\varepsilon}$, as $\varepsilon \rightarrow 0$. One then has three possible behaviors for the equations (1.5)-(1.6) in a limit $(N, g' , x, \tau )$ as $\varepsilon  \rightarrow $ 0, (in a subsequence). The discussion to follow here is formal in that we are not concerned  with the existence of such limits; this issue is discussed in detail in  [An4], as is the situation where $\tau = \tau_{\varepsilon}$ is not uniformly bounded in $B_\varepsilon$ as $\varepsilon \rightarrow 0$. (It turns out that the limits below have the same form even if $\{\tau_{\varepsilon}\}$ is unbounded).

\medskip

{\bf Case(i).}
 $\varepsilon / \rho^{2} \rightarrow $ 0.

 In this case, the equations (1.7)-(1.8) in the limit $\varepsilon  \rightarrow $ 0 take the form
\begin{equation} \label{e1.8}
L^{*}(\tau ) = 0, \Delta\tau  = 0. 
\end{equation}
These are the static vacuum Einstein equations, c.f. [An2,3] or [EK].

{\bf Case(ii).}
 $\varepsilon / \rho^{2} \rightarrow  \alpha  > $ 0.

 In this case, the limit equations take the form
\begin{equation} \label{e1.9}
\alpha\nabla{\cal R}^{2} + L^{*}(\tau ) = 0, 
\end{equation}
\begin{equation} \label{e1.10}
\Delta\tau  = -\frac{\alpha}{4}|r|^{2}. 
\end{equation}
Formally, these are the equations for a critical metric $g' $ of ${\cal R}^{2}$ subject to the constraint that $s =$ 0. However there are no compact scalar-flat perturbations of $g' $ since the limit $(N, g' )$ is non-compact and thus one must impose certain boundary conditions on the comparison metrics. For example, if the limit $(N, g' )$ is complete and asymptotically flat, (1.9)-(1.10) are the equations for a critical point of ${\cal R}^{2}$ among all scalar-flat and asymptotically flat metrics with a given mass $m$.

{\bf Case(iii).}
 $\varepsilon / \rho^{2} \rightarrow  \infty .$

 In this case, renormalize the equations (1.5)-(1.6) by dividing by $\varepsilon /\rho^{2}.$ Since $\tau $ is bounded, $(\rho^{2}/\varepsilon )\tau  \rightarrow $ 0, and one obtains in the limit
\begin{equation} \label{e1.11}
\nabla{\cal R}^{2} = 0, 
\end{equation}
\begin{equation} \label{e1.12}
|r|^{2} = 0, 
\end{equation}
so the limit metric is flat.

 These three cases may be summarized by the equations
\begin{equation} \label{e1.13}
\alpha\nabla{\cal R}^{2} + L^{*}(\tau ) = 0, 
\end{equation}
\begin{equation} \label{e1.14}
\Delta\tau  = -\frac{\alpha}{4}|r|^{2}, 
\end{equation}
where $\alpha  =$ 0 corresponds to Case (i), 0 $<  \alpha  <  \infty $ corresponds to Case (ii) and $\alpha  = \infty $ corresponds to (the here trivial) Case (iii). 

 Essentially the same discussion is valid for the scale-invariant functional
\begin{equation} \label{e1.15}
J_{\varepsilon} = (\varepsilon v^{1/3}{\cal R}^{2} -  v^{2/3}\cdot  s)|_{{\cal C}}, 
\end{equation}
where ${\cal C} $ is the space of Yamabe metrics on $M$. The existence and general properties of minimizers of $J_{\varepsilon}$ again are discussed in [An1,\S 8II]. By the same considerations, one obtains as above limit equations of the form (1.13)-(1.14), with $\tau $ replaced by the potential function $- (1+h)$ from [An1,\S 8II].

 Next consider briefly the scale-invariant functional
\begin{equation} \label{e1.16}
I_{\varepsilon}' = \varepsilon v^{1/3}\int|r|^{2} + \bigl( v^{1/3}\int (s^-)^{2}\bigr )^{1/2}, 
\end{equation}
on the space of metrics on $M$ as above, where $s^-$ = min$(s, 0)$. This functional, (essentially), is the main focus of [An4], and we refer there for a complete discussion. c.f. also \S 5.2. The Euler-Lagrange equations of $I_{\varepsilon}^-$ are formally the same as (1.2)-(1.3) with $\tau^-$ = min$(\tau ,0)$ in place of $\tau $ and $\sigma $ replaced by $(v^{1/3}\int (s^-)^{2})^{1/2}.$ Let now $g_{\varepsilon}$ be a minimizer of $I_{\varepsilon}^-$ on $M$. Note that here $\tau^-$ is automatically bounded above, as is $\int (s^-)^{2}$ as $\varepsilon  \rightarrow $ 0. However, there is no longer an apriori bound on the $L^{2}$ norm of $s$, i.e. it may well happen that $\int s^{2}(g_{\varepsilon}) \rightarrow  \infty $ as $\varepsilon  \rightarrow $ 0.

 Formally taking a blow-up limit $(N, g' , x, \tau^-)$ as $\varepsilon  \rightarrow $ 0 as above leads to the analogue of the equations (1.13)-(1.14), i.e. to the limit equations
\begin{equation} \label{e1.17}
\alpha\nabla{\cal R}^{2} + L^{*}(\tau^-) = 0, 
\end{equation}
\begin{equation} \label{e1.18}
\Delta (\tau^- + {\tfrac{3}{4}}\alpha s) = -{\tfrac{1}{4}}\alpha |r|^{2}, 
\end{equation}
where as before we assume that $\tau^-$ = lim $\tau_{\varepsilon}^-$ is bounded below. These equations correspond to (0.4)-(0.5) with $\omega  = \tau^-,$ but with $s \geq $ 0 in place of $s =$ 0, corresponding to the fact that in blow-up limits, now only $\int (s^-)^{2} \rightarrow $ 0 while previously $\int s^{2} \rightarrow $ 0.

 While in the region $N^- = {\tau^- <  0}$, the equations (1.17)-(1.18) have the same form as the Cases (i)-(iii) above, in the region $N^{+} = \{s > $ 0\}, (so that $\tau^- =$ 0), these equations take the form
\begin{equation} \label{e1.19}
\nabla{\cal R}^{2} = 0, 
\end{equation}
\begin{equation} \label{e1.20}
\Delta s = -\tfrac{1}{3}|r|^{2}, 
\end{equation}
i.e. the ${\cal R}^{2}$ equations (0.2)-(0.3); here we have divided by $\alpha .$ 

 The junction $\Sigma  = \partial\{s =$ 0\} between the two regions $N^-$ and $N^{+}$ above is studied in [An4]. This junction may be compared with junction conditions common in general relativity, where vacuum regions of space-(time) are joined to regions containing a non-vanishing matter distribution, c.f. [W, Ch.6.2] or [MTW, Ch.21.13, 23].

 This concludes the brief discussion on the origin of the equations in \S 0. The remainder of the paper is concerned with properties of their solutions.

\section{Non-Existence of ${\cal R}^{2}$ Solutions.}
\setcounter{equation}{0}

 In this section, we prove Theorem 0.1, i.e. there are no non-trivial complete ${\cal R}^{2}$ solutions with non-negative scalar curvature. The proof will proceed in several steps following in broad outline the proofs of [An1, Thms. 6.1,6.2].

  We begin with some preliminary material. Let $r_{h}(x)$ and $\rho (x)$ denote the $L^{2,2}$ harmonic radius and $L^{2}$ curvature radius of $(N, g)$ at $x$, c.f. [An1, Def.3.2] for the exact definition. Roughly speaking, $r_{h}(x)$ is the largest radius of the geodesic ball at $x$ on which there exists a harmonic chart for $g$ in which the metric differs from the flat metric by a fixed small amount, say $c_{o},$ in the $L^{2,2}$ norm. Similarly, $\rho (x)$ is the largest radius on which the $L^{2}$ average of the curvature is bounded by $c_{o}\cdot \rho (x)^{-2},$ for some fixed but small constant $c_{o} > $ 0. From the definition, 
\begin{equation} \label{e2.1}
\rho (y) \geq  dist(y, \partial B_{x}(\rho (x))), 
\end{equation}
for all $y\in B_{x}(\rho (x)),$ and similarly for $r_{h}(x).$ The point $x$ is called (strongly) $(\rho ,d)$ buffered if $\rho (y) \geq  d\cdot \rho (x),$ for all $y\in\partial B_{x}(\rho (x)),$ c.f. [An3, Def.3.7] and also [An1, \S 5]. This condition insures that there is a definite amount of curvature in $L^{2}$ away from the boundary in $B_{x}(\rho (x)).$

 As shown in [An1,\S 4], the equations (0.2)-(0.3) form an elliptic system and hence satisfy elliptic regularity estimates. Thus within the $L^{2,2}$ harmonic radius $r_{h},$ one actually has $C^{\infty}$ bounds of the solution metric, and hence its curvature, away from the boundary; the bounds depend only on the size of $r_{h}^{-1}.$ Observe also that the ${\cal R}^{2}$ equations are scale-invariant.

 Given this regularity derived from (0.2), as mentioned in \S 0, the full equation (0.2) itself is not otherwise needed for the proof of Theorem 0.1; only the trace equation (0.3) is used from here on.

\medskip

  For completeness, we recall some results from [An1, \S 2, \S 3] concerning convergence and collapse of Riemannian manifolds with uniform lower bounds on the $L^2$ curvature radius. Thus, suppose ${(B_{i}(R), g_i, x_i)}$ is a sequence of geodesic $R$-balls in complete non-compact Riemannian manifolds $(N_i, g_i)$, centered at base points $x_i \in N_i$. Suppose that $\rho_{i}(y_{i}) \geq \rho_{o}$, for some $\rho_{o} > 0$, for all $y_i \in B_{x_{i}}(R)$. The sequence is said to be {\it non-collapsing} if there is a constant $\nu_o > 0$ such that vol$B_{x_{i}}(1) \geq \nu_o$, for all $i$, or equivalently, the volume radius of $x_i$ is uniformly bounded below. In this case, it follows that a subsequence of ${(B_{i}(R), g_i, x_i)}$ converges in the weak $L^{2,2}$ topology to a limit ${(B_{\infty}(R), g_{\infty}, x_{\infty})}$. The convergence is uniform on compact subsets, and the limit is a manifold with $L^{2,2}$ Riemannian metric $g_{\infty}$, with base point $x_{\infty}$ = lim $x_i$. In case the sequence above is a sequence of ${\cal R}^2$ solutions, the convergence above is in the $C^{\infty}$ topology, by the regularity results above.

  The sequence as above is {\it collapsing} if vol$B_{x_{i}}(1) \rightarrow 0$, as $i \rightarrow \infty$. In this case, it follows that vol$B_{y_{i}}(1) \rightarrow 0$, for all $y_i \in B_{x_{i}}(R)$. Further, for any $\delta > 0$, and $i$ sufficiently large, there are domains $U_{i} = U_{i}(\delta)$, with $B_{x_{i}}(R - 2\delta) \subset U_i \subset B_{x_{i}}(R - \delta)$ such that, topologically, $U_i$ is a graph manifold. Thus, $U_i$ admits an F-structure, and the $g_i$-diameter of the fibers, (circles or tori), converges to 0 as $i \rightarrow \infty$. Now in case the sequence satisfies regularity estimates as stated above for ${\cal R}^2$ solutions, the curvature is uniformly bounded in $L^{\infty}$ on $U_i$. Hence, it follows from results of Cheeger-Gromov, Fukaya and Rong, (c.f. [An1, Thm.2.10]), that $U_i$ is topologically either a Seifert fibered space or a torus bundle over an interval. Further, the inclusion map of any fiber induces an injection into the fundamental group $\pi_{1}(U_i)$, and the fibers represent (homotopically) the collection of all very short essential loops in $U_i$. Hence, there is an infinite ${\Bbb Z}$ or ${\Bbb Z} \oplus {\Bbb Z}$ cover $(\widetilde U_i, g_i, \widetilde x_i)$ of $(U_i, g_i, x_i)$, ($\widetilde x_i$ a lift of $x_i$), obtained by unwrapping the fibers, which is a non-collapsing sequence. For if the lifted sequence of covers collapsed, by the same arguments again, $\widetilde U_i$, and hence $U_i$, must contain essential short loops; all of these however have already been unwrapped in $\widetilde U_i$. Alternately, one may pass to sufficiently large finite covers $\bar U_i$ of $U_i$ to unwrap the collapse, in place of the infinite covers. In this collapse situation, the limit metrics $(\widetilde U_{\infty}, g_{\infty}, \widetilde x_{\infty})$, (resp. $(\bar U_{\infty}, g_{\infty}, \bar x_{\infty})$) have free isometric ${\Bbb R}$, (resp. $S^1$), actions.

  These results on the behavior of non-collapsing and collapsing sequences will be used frequently below.

\medskip

 We now begin with the proof of Theorem 0.1 itself. The following Lemma shows that one may assume without loss of generality that a complete ${\cal R}^{2}$ solution has uniformly bounded curvature.
\bbgin{lemma} \label{l 2.1.}
  Let $(N, g)$ be a complete non-flat ${\cal R}^{2}$ solution. Then there exists another complete non-flat ${\cal R}^{2}$ solution $(N' , g' ),$ obtained as a geometric limit at infinity of $(N, g)$, which has uniformly bounded curvature, i.e.
\begin{equation} \label{e2.2}
|r|_{g'} \leq  1. 
\end{equation}

\end{lemma}
{\bf Proof:}
 We may assume that $(N, g)$ itself has unbounded curvature, for otherwise there is nothing to prove. It follows from the $C^{\infty}$ regularity of solutions mentioned above that the curvature $|r|$ is unbounded on a sequence $\{x_{i}\}$ in $(N, g)$ if and only if
\begin{equation} \label{e2.3}
\rho (x_{i}) \rightarrow  0. 
\end{equation}
For such a sequence, let $B_{i} = B_{x_{i}}(1)$ and let $d_{i}(x) = dist(x_{i}, \partial B_{i}).$ Consider the scale-invariant ratio $\rho (x)/d_{i}(x),$ for $x\in B_{i},$ and choose points $y_{i}\in B_{i}$ realizing the minimum value of $\rho /d_{i}$ on $B_{i}.$ Since $\rho /d_{i}$ is infinite on $\partial B_{i}, y_{i}$ is in the interior of $B_{i}.$ By (2.3), we have 
$$\rho (y_{i})/d_{i}(y_{i}) \rightarrow  0, $$
and so in particular $\rho (y_{i}) \rightarrow $ 0. Now consider the sequence $(B_{i}, g_{i}, y_{i}),$ where $g_{i} = \rho (y_{i})^{-2}\cdot  g.$ By construction, $\rho_{i}(y_{i}) =$ 1, where $\rho_{i}$ is the $L^{2}$ curvature radius w.r.t. $g_{i}$ and $\delta_{i}(y_{i}) = dist_{g_{i}}(y_{i}, \partial B_{i}) \rightarrow  \infty .$ Further, by the minimality property of $y_{i},$
\begin{equation} \label{e2.4}
\rho_{i}(x) \geq  \rho_{i}(y_{i})\cdot \frac{\delta_{i}(x)}{\delta_{i}(y_{i})} = \frac{\delta_{i}(x)}{\delta_{i}(y_{i})} . 
\end{equation}
It follows that $\rho_{i}(x) \geq  \frac{1}{2},$ at all points $x$ of uniformly bounded $g_{i}$-distance to $y_{i},$ (for $i$ sufficiently large, depending on $dist_{g_{i}}(x, y_{i})).$ 

 Consider then the pointed sequence $(B_{i}, g_{i}, y_{i}).$ If this sequence, (or a subsequence), is not collapsing at $y_{i},$ then the discussion above, applied to ${(B_{y_{i}}(R_j), g_i, y_i)}$, with $R_j \rightarrow \infty$, implies that a diagonal subsequence converges smoothly to a limit $(N' , g' , y)$, $y =$ lim $y_{i}.$ The limit is a complete ${\cal R}^{2}$ solution, (since $\delta_{i}(y_{i}) \rightarrow  \infty )$ satisfying $\rho  \geq  \frac{1}{2}$ everywhere, and $\rho (y) = 1$, since $\rho$ is continuous under smooth convergence to limits, (c.f. [An1, Thm. 3.5]). Hence the limit is not flat. By the regularity estimates above, the curvature $|r|$ is pointwise bounded above. A further bounded rescaling then gives (2.2).

 On the other hand, suppose this sequence is collapsing at $y_{i}.$ Then from the discussion preceding Lemma 2.1, it is collapsing everywhere within $g_{i}$-bounded distance to $x_{i}$ along a sequence of injective F-structures. Hence one may pass to suitable covers $\widetilde U_{i}$ of $U_i$ with $B_{x_{i}}(R_i - 1) \subset U_i \subset B_{x_{i}}(R)$, for some sequence $R_{i} \rightarrow  \infty $ as $i \rightarrow  \infty .$ This sequence is not collapsing and thus one may apply the reasoning above to again obtain a limit complete non-flat ${\cal R}^{2}$ solution satisfying (2.2), (which in addition has a free isometric ${\Bbb R}$-action.

{\endproof}

 Let $v(r)$ = vol$B_{x_{o}}(r),$ where $B_{x_{o}}(r)$ is the geodesic $r$-ball about a fixed point $x_{o}$ in $(N, g)$. Let $J^{2}$ be the Jacobian of the exponential map exp: $T_{x_{o}}N \rightarrow  N$, so that
$$v(r) = \int_{S_{o}}\int_{0}^{r}J^{2}(s,\theta )dsd\theta , $$
where $S_{o}$ is the unit sphere in $T_{x_{o}}N.$ Thus,
$$v' (r) = \int_{S_{o}}J^{2}(r,\theta )d\theta , $$
is the area of the geodesic sphere $S_{x_{o}}(r).$

 The next result proves Theorem 0.1 under reasonably weak conditions, and will also be needed for the proof in general.
\bbgin{proposition} \label{p 2.2.}
  Let (N, g) be a complete ${\cal R}^{2}$ solution on a 3-manifold N, with bounded curvature and $s \geq $ 0. Suppose there are constants $\varepsilon  > $ 0 and $c <  \infty $ such that
\begin{equation} \label{e2.5}
v(r) \leq  c\cdot  r^{4-\varepsilon}, 
\end{equation}
for all $r \geq $ 1. Then (N, g) is flat.
\end{proposition}
{\bf Proof:}
 Let $t(x)$= dist$(x, x_{o})$ be the distance function from $x_{o}\in N$ and let $\eta  = \eta (t)$ be a non-negative cutoff function, of compact support to be determined below, but initially satisfying $\eta' (t) \leq $ 0. Multiply (0.3) by $\eta^{4}$ and apply the divergence theorem, (this is applicable since $\eta (t)$ is a Lipschitz function on $N$), to obtain
$$\int\eta^{4}|r|^{2} = 3\int<\nabla s, \nabla\eta^{4}> . $$
Now one cannot immediately apply the divergence theorem again, since $t$ and hence $\eta $ is singular at the cut locus $C$ of $x_{o}.$ Let $U_{\delta}$ be the $\delta$-tubular neighborhood of $C$ in $N$. Then applying the divergence theorem on $N \setminus U_{\delta}$ gives
$$\int_{N \setminus {U_{\delta}}}<\nabla s, \nabla\eta^{4}>  = -\int_{N \setminus {U_{\delta}}}s\Delta\eta^{4} + \int_{\partial (N \setminus {U_{\delta}})}s<\nabla\eta^{4}, \nu> , $$
where $\nu $ is the unit outward normal. Since $<\nu , \nabla t>  > $ 0 on $\partial (N \setminus {U_{\delta}})$ and $\eta'  \leq $ 0, the hypothesis $s \geq $ 0 implies that the boundary term is non-positive. Hence,
$$\int_{N \setminus U_{\delta}} \eta^{4}|r|^{2} \leq  - 3\int_{N \setminus U_{\delta}} s\Delta\eta^{4}.$$
We have $\Delta\eta^{4} = 4\eta^{3}\Delta\eta  + 12\eta^{2}|\nabla \eta|^{2} \geq  4\eta^{3}\Delta\eta ,$ so that again since $s \geq $ 0,
$$\int_{N \setminus U_{\delta}} \eta^{4}|r|^{2} \leq  - 12\int_{N \setminus U_{\delta}} s\eta^{3}\Delta\eta . $$
Further $\Delta\eta  = \eta'\Delta t + \eta'' ,$ so that
$$\int_{N \setminus U_{\delta}} \eta^{4}|r|^{2} \leq  - 12\int_{N \setminus U_{\delta}} s\eta^{3}\eta' \Delta t + s\eta^{3}\eta''. $$
It is standard, c.f. [P,9.1.1], that off $C$,
$$\Delta t = H = 2\frac{J'}{J}, $$
where $H$ is the mean curvature of $S_{x_{o}}(r)$ and $J = (J^{2})^{1/2}.$ Further, since the curvature of $(N, g)$ is bounded, standard comparison geometry, (c.f. [Ge] for example), implies that there is a constant $C < \infty$ such that
\begin{equation} \label{e2.6}
H(x) \leq C,
\end{equation}
for all $x$ outside $B_{x_{o}}(1) \subset N$. (Of course there is no such lower bound for $H$). Hence, since $s \geq 0$ and $\eta' \leq 0$, it follows that
$$\int_{N \setminus U_{\delta}} \eta^{4}|r|^{2} \leq  - 12\int_{N \setminus U_{\delta}} s\eta^{3}\eta' H^{+} + s\eta^{3}\eta'',$$
where $H^{+}$ = max$(H, 0)$. The integrand $-s\eta'H^{+}$ is positive and bounded. Hence, since the cutlocus $C$ is of measure 0, we may let $\delta \rightarrow 0$ and obtain
$$\int_{N} \eta^{4}|r|^{2} \leq  - 12\int_{N} s\eta^{3}\eta' H^{+} + s\eta^{3}\eta''.$$

  Now fix any $R <  \infty $ and choose $\eta  = \eta (t)$ so that $\eta  \equiv $ 1 on $B_{x_{o}}(R), \eta  \equiv $ 0 on $N \setminus B_{x_{o}}(2R), \eta'  \leq $ 0, and $|\eta'| \leq  c/R, |\eta''| \leq  c/R^{2}.$ Using the H\"older and Cauchy inequalities, we obtain
\begin{equation} \label{e2.7}
\int\eta^{4}|r|^{2} \leq  \mu\int\eta^{4}s^{2} + \mu^{-1}\int\eta^{2}(\eta' )^{2}(H^{+})^{2} + \mu^{-1}\int\eta^{2}(\eta'' )^{2},
\end{equation}
for any $\mu  > $ 0 small. Since $|r|^{2} \geq  s^{2}/3,$ by choosing $\mu $ sufficiently small the first term on the right in (2.7) may be absorbed into the left. Thus we have on $B(R)$= $B_{x_{o}}(R),$ for suitable constants $c_{i}$ independent of $R$,
\begin{equation} \label{e2.8}
\int_{B(R)}|r|^{2} \leq  c_{1}\int_{B(2R)}(R^{-2}(H^{+})^{2} + R^{-4}) \leq  c_{2}R^{-2}\int_{B(2R)}(H^{+})^{2} + c_{3}R^{-\varepsilon}, 
\end{equation}
where the last inequality uses (2.5). 

  We now claim that there is a constant $K < \infty$, (depending on the geometry of $(N, g)$), such that
\begin{equation} \label{e2.9}
\Delta t(x) \cdot \rho(x) \leq K,
\end{equation} 
for all $x \in N$ with $t(x) \geq 10$, with $x \notin C$. We will assume (2.9) for the moment and complete the proof of the result; following this, we prove (2.9).

  Thus, substituting (2.9) in (2.8), and using the definition of $\rho$, we obtain 
$$\int_{B(R)}|r|^{2}  \leq  c_{4}R^{-2}\int_{B(2R)}|r| + c_{3}R^{-\varepsilon}.$$
Applying the Cauchy inequality to the $|r|$ integral then gives
$$\int_{B(R)}|r|^{2} \leq  c_{4}R^{-2}\big (\int_{B(2R)}|r|^{2} \bigl )^{1/2}vol B(2R)^{1/2} + c_{3}R^{-\varepsilon}.$$
Now from the volume estimate (2.5) and the uniform bound on $|r|$, there exists a sequence $R_i \rightarrow \infty$ and a constant $C < \infty$ such that
$$\int_{B(2R_i)}|r|^{2} \leq C\int_{B(R_i)}|r|^{2}.$$ 
Hence, setting $R = R_i$ and combining these estimates gives
$$\int_{B(R_i)}|r|^{2}  \leq c_{5}\frac{vol B(2R_i)}{R_{i}^{4}}.$$
Taking the limit as $i \rightarrow \infty$ and using (2.5), it follows that $(N, g)$ is flat, as required.

  Thus, it remains to establish (2.9). We prove (2.9) by contradiction. Thus, suppose there is a sequence $\{x_i\} \in N \setminus C$ such that
\begin{equation} \label{e2.10} 
\Delta t(x_i) \cdot \rho(x_i) \rightarrow \infty,
\end{equation}
as $i \rightarrow \infty$. Note that (2.10) is scale invariant and that necessarily $t(x_i) \rightarrow \infty$. In fact, by (2.6), (2.10) implies that $\rho(x_i) \rightarrow \infty$ also. Note that
\begin{equation} \label{e2.11}
\rho(y) \leq 2t(y),
\end{equation}
for any $y$ such that $t(y)$ is sufficiently large, since $(N, g)$ is assumed not flat.

  We rescale the manifold $(N, g)$ at $x_i$ by setting $g_i = \lambda_{i}^{2} \cdot g$, where $\lambda_{i} = \Delta t(x_i)$. Thus, w.r.t. $g_i$, we have $\Delta_{g_{i}}t_{i}(x_{i}) = 1$, where $t_{i}(y)$ = $dist_{g_{i}}(y, x_o) = \lambda_i t(y)$. By the scale invariance of (2.10), it follows that
\begin{equation} \label{e2.12}
\rho_i(x_i) \rightarrow \infty ,
\end{equation}
where $\rho_i = \lambda_i \cdot \rho$ is the $L^2$ curvature radius w.r.t. $g_i$. By (2.11), this implies that $t_i(x_i) \rightarrow \infty$, so that the base point $x_o$ diverges to infinity in the $\{x_i\}$ based sequence $(N, g_i, x_i)$. Hence renormalize $t_i$ by setting $\beta_i(y) = t_i(y) - dist_{g_{i}}(y, x_o)$, as in the construction of Busemann functions.

  Thus, we have a sequence of ${\cal R}^2$ solutions $(N, g_i, x_i)$ based at $\{x_i\}$. From the discussion preceding Lemma 2.1, it follows that a subsequence converges smoothly to an ${\cal R}^2$ limit metric $(N_{\infty}, g_{\infty}, x_{\infty})$, passing to suitable covers as described in the proof of Lemma 2.1 in the case of collapse. By (2.12), it follows that
$$N_{\infty} = {\Bbb R}^3, $$
(or a quotient of ${\Bbb R}^3$), and $g_{\infty}$ is the complete flat metric. Now the smooth convergence also gives
$$\Delta_{g_{\infty}}\beta(x_{\infty}) = 1, $$
where $\beta$, the limit of $\beta_i$, is a Busemann function on a complete flat manifold. Hence $\beta$ is a linear coordinate function. This of course implies $\Delta_{g_{\infty}}\beta(x_{\infty}) = 0$, giving a contradiction. This contradiction then establishes (2.9).

{\endproof}

 We remark that this result mainly requires the hypothesis $s \geq $ 0 because of possible difficulties at the cut locus. There are other hypotheses that allow one to overcome this problem. For instance if $(N, g)$ is complete as above and (2.5) holds, (but without any assumption on $s$), and if there is a smooth approximation $\Roof{t}{\widetilde}$ to the distance function $t$ such that $|\Delta\Roof{t}{\widetilde}| \leq  c/\Roof{t}{\widetilde},$ (for example if $|r| \leq  c/t^{2}$ for some $c <  \infty ),$ then $(N, g)$ is flat. The proof is the same as above, (in fact even simpler in this situation).

\medskip
  Next we need the following simple result, which allows one to control the full curvature in terms of the scalar curvature. This result is essentially equivalent to [An1, Lemma 5.1].
\bbgin{lemma} \label{l 2.3.}
  Let $g$ be an ${\cal R}^{2}$ solution, defined in a geodesic ball $B = B_{x}(1),$ with $r_{h}(x) =$ 1. Then for any small $\mu  > $ 0, there is a constant $c_{1} = c_{1}(\mu )$ such that
\begin{equation} \label{e2.13}
|r|^{2}(y) \leq  c_{1}\cdot ||s||_{L^{2}(B)} , 
\end{equation}
for all $y\in B(1-\mu ) = B_{x}(1-\mu ).$ In particular, if $||s||_{L^{2}(B)}$ is sufficiently small, then $g$ is almost flat, i.e. has almost 0 curvature, in $B(1-\mu ).$

 Further, if $s \geq $ 0 in $B(1)$, then there is a constant $c_{2} = c_{2}(\mu )$ such that
\begin{equation} \label{e2.14}
||s||_{L^{2}(B(1-\mu ))} \leq  c_{2}s(x). 
\end{equation}
\end{lemma}

{\bf Proof:}
 Let $\eta $ be a non-negative cutoff function satisfying $\eta  \equiv $ 1 on $B(1-\frac{\mu}{2}), \eta  \equiv $ 0 on $A(1-\frac{\mu}{4},1),$ and $|\nabla \eta| \leq  c/\mu .$ Pair the trace equation (0.3) with $\eta^{2}$ to obtain
$$\int_{B}\eta^{2}|r|^{2} = - 3\int_{B}s\Delta\eta  \leq  c\cdot  (\int_{B}s^{2})^{1/2}(\int_{B}(\Delta\eta )^{2})^{1/2}. $$
Since $r_{h}(x) =$ 1, $\eta $ may be chosen so that the $L^{2}$ norm of $\Delta\eta $ is bounded in terms of $\mu $ only. It follows that
$$\int_{B(1-\frac{\mu}{2})}|r|^{2} \leq  c(\mu )||s||_{L^{2}(B)}. $$
One obtains then an $L^{\infty},$ (and in fact $C^{k,\alpha}),$ estimate for $|r|^{2}$ by elliptic regularity, as discussed preceding Lemma 2.1.

 For the second estimate (2.14), note that by (0.3), $s$ is a superharmonic function, assumed non-negative in $B(1)$. Since the metric $g$ is bounded in $L^{2,2}$ on $B(1)$, and hence bounded in $C^{1/2}$ by Sobolev embedding, the estimate (2.14) is an immediate consequence of the DeGiorgi-Nash-Moser estimates for non-negative supersolutions of divergence form elliptic equations, c.f. [GT,Thm.8.18].

{\endproof}

 The behavior of the scalar curvature, and thus of the full curvature, at infinity is the central focus of the remainder of the proof. For example, Lemma 2.3 leads easily to the following special case of Theorem 0.1.

\bbgin{lemma} \label{l 2.4.}
  Suppose (N, g) is a complete ${\cal R}^{2}$ solution satisfying 
\begin{equation} \label{e2.15}
limsup_{t\rightarrow\infty} \ t^{2}\cdot  s = 0, 
\end{equation}
where t(x) $=$ dist(x, $x_{o}).$ Then (N, g) is flat.
\end{lemma}
{\bf Proof:}
 We claim first that (2.15) implies that
\begin{equation} \label{e2.16}
liminf_{t\rightarrow\infty} \ \rho /t \geq  c_{o}, 
\end{equation}
for some constant $c_{o} > $ 0. For suppose (2.16) were not true. Then there is a sequence $\{x_{i}\}$ in $N$ with $t_{i} = t(x_{i}) \rightarrow  \infty ,$ such that $\rho (x_{i})/t(x_{i}) \rightarrow $ 0. We may choose $x_{i}$ so that it realizes approximately the minimal value of the ratio $\rho /t$ for $ \frac{1}{2}t_{i} \leq t \leq  2t_{i},$ as in the proof of Lemma 2.1. For example, choose $x_{i}$ so that it realizes the minimal value of the ratio
$$\rho (x)/dist(x, \partial A({\tfrac{1}{2}}t_{i}, 2t_{i})) $$
for $x\in A(\frac{1}{2}t_{i}, 2t_{i}).$ Such a choice of $x_{i}$ implies that $x_{i}$ is strongly $(\rho ,\frac{1}{2})$ buffered, i.e. $\forall y_{i}\in\partial B_{x_{i}}(\rho (x_{i})),$
$$\rho (y_{i}) \geq  \tfrac{1}{2}\rho (x_{i}), $$
c.f. the beginning of \S 2 and compare with (2.4).

 Now rescale the metric $g$ by the $L^{2}$ curvature radius $\rho $ at $x_{i},$ i.e. set $g_{i} = \rho (x_{i})^{-2}\cdot  g.$ Thus $\rho_{i}(x_{i}) =$ 1, where $\rho_{i} = \rho (g_{i}).$ As in the proof of Lemma 2.1, if the ball $(B_{i}, g_{i}), B_{i} = (B_{x_{i}}(\frac{11}{8}), g_{i})$ is sufficiently collapsed, pass to sufficiently large covers of this ball to unwrap the collapse, as discussed preceding Lemma 2.1. We assume this is done, and do not change the notation for the collapse case.

 Since we are assuming that $\rho (x_{i}) <<  t(x_{i}),$ by (2.15) and scaling properties, it follows that $s_{i},$ the scalar curvature of $g_{i},$ satisfies
\begin{equation} \label{e2.17}
s_{i} \rightarrow  0, 
\end{equation}
uniformly on $(B_{i}(\frac{5}{4}), g_{i})$. By Lemma 2.3, we obtain
\begin{equation} \label{e2.18}
|r_{i}| \rightarrow  0, 
\end{equation}
uniformly on $(B_{x_{i}}(\frac{9}{8}), g_i)$. However, since $\rho_{i}(x_{i}) =$ 1, and vol$B_{x{_i}}(1) > \nu_o > 0$, the ball $(B_{x_{i}}(\frac{9}{8}), g_i)$ has a definite amount of curvature in $L^{2}.$ This contradiction gives (2.16).

 Now apply the same reasoning to any sequence $\{y_{i}\}$ in $N$, with $t(y_{i}) \rightarrow  \infty ,$ but with respect to the blow-down metrics $g_{i} = t(y_{i})^{-2}\cdot  g,$ so that by (2.16), $\rho_{i}(y_{i}) \geq  c > $ 0. Since (2.17) remains valid, apply Lemma 2.3 again to obtain (2.18) on a neighborhood of fixed $g_{i}$-radius about $y_{i}.$ The estimate (2.18) applied to the original (unscaled) metric $g$ means that
\begin{equation} \label{e2.19}
limsup_{t\rightarrow\infty} \ t^{2}\cdot |r| = 0, 
\end{equation}
improving the estimate (2.15).

 Now standard comparison estimates on the Ricatti equation $H' + \frac{1}{2}H^2 \leq |r|$, (c.f. [P, Ch.9] and the proof of Prop. 2.2), shows that (2.19) implies that the volume growth of $(N, g)$ satisfies
$$v(t) \leq  c\cdot  t^{3+\varepsilon}, $$
for any given $\varepsilon  > $ 0 with $c = c(\varepsilon ) <  \infty .$ The maximum principle applied to the trace equation (0.3), together with (2.15) implies that $s > 0$ everywhere. Thus, the result follows from Proposition 2.2.

{\endproof}

 The proof of Theorem 0.1 now splits into two cases, following the general situation in [An1, Thms. 6.1, 6.2] respectively. The first case below can basically be viewed as a local and quantitative version of Lemma 2.4. The result roughly states that if a complete ${\cal R}^{2}$ solution with $s \geq $ 0 is weakly asymptotically flat in some direction, then it is flat. 
\bbgin{theorem} \label{t 2.5.}
  Let (N, g) be a complete ${\cal R}^{2}$ solution with non-negative scalar curvature. Suppose there exists a sequence $x_{i}$ in (N, g) such that
\begin{equation} \label{e2.20}
\rho^{2}(x_{i})\cdot  s(x_{i}) \rightarrow  0 \ \ {\rm as} \ i \rightarrow  \infty . 
\end{equation}
Then (N, g) is flat.
\end{theorem}
{\bf Proof:}
 The proof follows closely the ideas in the proof of [An1, Thms. 5.4 and 6.1]. Throughout the proof below, we let $\rho $ denote the $L^{4}$ curvature radius as opposed to the $L^{2}$ curvature radius. As noted in [An1,(5.6)], the $L^{2}$ and $L^{4}$ curvature radii are uniformly equivalent to each other on $(N, g)$, since as discussed preceding Lemma 2.1, the metric satisfies an elliptic system, and regularity estimates for such equations give $L^{4}$ bounds in terms of $L^{2}$ bounds. In particular, Lemma 2.3 holds with the $L^{4}$ curvature radius in place of the $L^{2}$ radius. Further, as discussed preceding Lemma 2.1, if a ball $B(\rho )\subset (N, g)$ is sufficiently collapsed, we will always assume below that the collapse is unwrapped by passing to the universal cover. Thus, the $L^{4}$ curvature radius and $L^{2,4}$ harmonic radius are uniformly equivalent to each other, c.f. [An1, (3.8)-(3.9)].

 Let $\{x_{i}\}$ be a sequence satisfying (2.20). As in the proof of Lemma 2.4, (c.f. (2.17)ff), a subsequence of the rescaled metrics $g_{i} = \rho (x_{i})^{-2}\cdot  g$ converges to a flat metric on uniformly compact subsets of $(B_{x_{i}}(1), g_{i}),$ unwrapping to the universal cover in case of collapse. Thus, the $(L^{4})$ curvature radius $\rho_{i}(y_{i})$ w.r.t. $g_{i}$ necessarily satisfies $\rho_{i}(y_{i}) \rightarrow $ 0, for some $y_{i}\in (\partial B_{x_{i}}(1), g_{i}),$ as $i \rightarrow  \infty .$

 Pick $i_{o}$ sufficiently large, so that $g_{i_{o}}$ is very close to the flat metric. We relabel by setting $g^{1} = g_{i_{o}}, q^{1} = x_{i_{o}}$ and $B = B^{1}= (B_{q^{1}}(1), g^{1}).$ For the moment, we work in the metric ball $(B^{1}, g^{1}).$ It follows that for any $\delta_{1}$ and $\delta_{2} > $ 0, we may choose $i_{o}$ such that
\begin{equation} \label{e2.21}
\rho (q) \leq  \delta_{1}\rho (q_{1}), 
\end{equation}
for some $q\in\partial B,$ and (from (2.20)),
\begin{equation} \label{e2.22}
s(q^{1}) \leq  \delta_{2}, 
\end{equation}
where $\rho $ and $s$ are taken w.r.t. $g_{1}.$ Here both $\delta_{1}$ and $\delta_{2}$ are assumed to be sufficiently small, (for reasons to follow), and further $\delta_{2}$ is assumed sufficiently small compared with $\delta_{1}$ but sufficiently large compared with $\delta_1{}^{2}.$ For simplicity, and to be concrete, we set
\begin{equation} \label{e2.23}
\delta_{2} = \delta_1{}^{3/2}, 
\end{equation}
and assume that $\delta_{1}$ is (sufficiently) small.

 We use the trace equation (0.3), i.e.
\begin{equation} \label{e2.24}
\Delta s  = -\tfrac{1}{3}|r|^{2}, 
\end{equation}
on $(B, g^{1})$ to analyse the behavior of $s$ in this scale near $\partial B;$ recall that the trace equation is scale invariant. From the Green representation formula, we have for $x\in B$
\begin{equation} \label{e2.25}
s(x) = \int_{\partial B}P(x,Q)d\mu_{Q} -  \int_{B}\Delta s(y)G(x,y)dV_{y}, 
\end{equation}
where $P$ is the Poisson kernel and $G$ is the positive Green's function for the Laplacian on $(B, g^{1}).$ As shown in [An1, Lemma 5.2], the Green's function is uniformly bounded in $L^{2}(B).$ The same holds for $\Delta s$ by (2.24), since the $L^{4}$ curvature radius of $g^{1}$ at $q^{1}$ is 1. Thus, the second term in (2.25) is uniformly bounded, i.e. there is a fixed constant $C_{o}$ such that
\begin{equation} \label{e2.26}
s(x) \leq  \int_{\partial B}P(x,Q)d\mu_{Q} + C_{o}. 
\end{equation}
The Radon measure $d\mu $ is a positive measure on $\partial B,$ since $s > $ 0 everywhere. Further, the total mass of $d\mu $ is at most $\delta_{2},$ by (2.22).

 By [An1, Lemma 5.3], the Poisson kernel $P(x,Q)$ satisfies
$$P(x,Q) \leq  c_{1}\cdot  t_{Q}(x)^{-2}, $$
where $t_{Q}(x)$ = dist$(x,Q)$ and $c_{1}$ is a fixed positive constant. Hence, for all $x\in B,$
\begin{equation} \label{e2.27}
s(x) \leq  c_{1}\cdot \delta_{2}\cdot  t^{-2}(x) + C_{o}, 
\end{equation}
where $t(x)$= dist$(x,\partial B).$ 

 Now suppose the estimate (2.27) can be improved in the sense that there exist points $q^{2}\in B^{1}$ s.t.
\begin{equation} \label{e2.28}
\rho^{1}(q^{2}) \leq  \delta_{1}, 
\end{equation}
and
\begin{equation} \label{e2.29}
s(q^{2}) \leq  \tfrac{1}{2}\delta_{2}\cdot  (\rho^{1}(q^{2}))^{-2} + C_{o}, 
\end{equation}
where $\rho^{1}(x) = \rho (x, g^{1}).$ Note that $\rho^{1}(x) \geq  t(x)$ by (2.1), so that the difference between (2.27) and (2.29) is only in the factors $c_{1}$ and $\frac{1}{2}.$

 We have
\begin{equation} \label{e2.30}
\delta_{2}\cdot  (\rho^{1}(q^{2}))^{-2} \geq  \delta_{2}\delta_{1}^{-2} >>  1, 
\end{equation}
where the last estimate follows from the assumption (2.23) on the relative sizes of $\delta_{1}$ and $\delta_{2}.$ Thus, the term $C_{o}$ in (2.29) is small compared with its partner in (2.29), and so
\begin{equation} \label{e2.31}
s(q^{2}) \leq  \delta_{2}\cdot  (\rho^{1}(q^{2}))^{-2}. 
\end{equation}

 We may then repeat the analysis above on the new scale $g^{2} = (\rho^{1}(q^{2}))^{-2}\cdot  g^{1}$ and the $g^{2}$ geodesic ball $B^{2} = B_{q^{2}}^{2}(1),$ so that $\rho^{2}(q^{2}) =$ 1. Observe that the product $\rho^{2}\cdot  s$ is scale-invariant, so that in the $g^{2}$ scale, $s$ is much smaller than $s$ in the $g^{1}$ scale. In the $g^{2}$ scale, (2.31) becomes the statement
\begin{equation} \label{e2.32}
s(q^{2}) \leq  \delta_{2}, 
\end{equation}
as in (2.22).

 We will show below in Lemma 2.6 that one may continue in this way indefinitely, i.e. as long as there exist points $q^{k}\in B^{k-1},$ with $\rho (q^{k}) \leq  \delta_{1}\rho (q^{k-1})$ as in (2.28), then there exist such points $q^{k}$ satisfying in addition
\begin{equation} \label{e2.33}
s(q^{k}) \leq  \delta_{2}, 
\end{equation}
where $s$ is the scalar curvature of $g^{k} = (\rho^{k-1}(q_{k}))^{-2}\cdot  g^{k-1},$ as in (2.22) or (2.32).

 On the one hand, we claim this sequence $\{q^{k}\}$ must terminate at some value $k_{o}.$ Namely, return to the original metric $(N, g)$. By construction, we have $\rho (q^{k}) \leq  \delta_1{}^{k}\rho (q^{1}).$ The value $\rho (q^{1})$ is some fixed number, (possibly very large), say $\rho (q^{1}) = C$, so that
\begin{equation} \label{e2.34}
\rho (q^{k}) \leq  C\cdot \delta_1{}^{k} \rightarrow  0, \ \ {\rm as} \ \ k \rightarrow  \infty . 
\end{equation}
However, observe that $dist_{g}(q^{k}, q^{1})$ is uniformly bounded, independent of $k$. Since $(N, g)$ is complete and smooth, $\rho $ cannot become arbitrarily small in compact sets of $N$. Hence (2.34) prevents $k$ from becoming arbitrarily large.

 On the other hand, if this sequence terminates at $q^{k}, k = k_{o},$ then necessarily
\begin{equation} \label{e2.35}
\rho (q) \geq  \delta_{1}\cdot \rho (q^{k}), 
\end{equation}
for all $q\in\partial B^{k}.$ However the construction gives $s(q^{k}) \leq  \delta_{2},$ where $s$ is the scalar curvature of $g_{k},$ with $\rho^{k}(q^{k}) =$ 1. This situation contradicts Lemma 2.3 if $\delta_{2}$ is chosen sufficiently small compared with $\delta_{1},$ i.e. in view of (2.23), $\delta_{1}$ is sufficiently small.

 It follows that the proof of Theorem 2.5 is completed by the following:
\bbgin{lemma} \label{l 2.6.}

 Let (N, g) be an ${\cal R}^{2}$ solution, $x\in N,$ and let $g$ be scaled so that $\rho (x) =$ 1, where $\rho $ is the $L^{4}$ curvature radius. Suppose that $\delta_{1}$ is sufficiently small, $\delta_{2} = \delta_1{}^{3/2},$
\begin{equation} \label{e2.36}
s(x) \leq  \delta_{2}, 
\end{equation}
and, for some $y_{o}\in B_{x}(1),$
\begin{equation} \label{e2.37}
\rho (y_{o}) \leq  \delta_{1}\cdot \rho (x) = \delta_{1}. 
\end{equation}
  Then there exists an absolute constant $K <  \infty $ and a point $y_{1}\in B_{x}(1),$ with $\rho (y_{1}) \leq  K\cdot \rho (y_{o}),$ such that
\begin{equation} \label{e2.38}
s(y_{1}) \leq  \delta_{2}\cdot  (\rho (y_{1}))^{-2}. 
\end{equation}

\end{lemma}
{\bf Proof:}
  By (2.27), we have
\begin{equation} \label{e2.39}
s(y) \leq  c_{1}\cdot \delta_{2}\cdot  t^{-2}(y) + C_{o}, 
\end{equation}
for all $y\in B.$ If there is a $y$ in $B_{y_{o}}(2\rho (y_{o})) \cap B_x(1)$ such that (2.38) holds at $y$, then we are done, so suppose there is no such $y$. Consider the collection $\beta$ of points $z\in B$ for which an opposite inequality to (2.38) holds, i.e.
\begin{equation} \label{e2.40}
s(z) \geq  \tfrac{1}{10}\delta_{2}\cdot  t^{-2}(z), 
\end{equation}
for $t(z)$ very small; (the factor $\frac{1}{10}$ may be replaced by any other small positive constant). From (2.26), this implies that $s$ resembles a multiple of the Poisson kernel near $z$. More precisely, suppose (2.40) holds for all $z$ within a small ball $B_{z_{o}}(\nu ),$ for some $z_{o}\in\partial B.$ Given a Borel set $E \subset  \partial B,$ let $m(E)$ denote the mass of the measure $d\mu $ from (2.26), of total mass at most $\delta_{2}.$ Then the estimate (2.40) implies
\begin{equation} \label{e2.41}
m(B_{z_{o}}(\nu )) \geq  \delta_{2}\cdot \varepsilon_{o}, 
\end{equation}
where $\nu $ may be made arbitrarily small if (2.40) holds for $z\in\beta$ and $t(z)$ is sufficiently small. The constant $\varepsilon_{o}$ depends only on the choice of $\frac{1}{10}$ in (2.40). Thus, part of $d\mu $ is weakly close to a multiple of the Dirac measure at some point $z_{o}\in\partial B$ near $\beta$; (the Dirac measure at $z_{o}$ generates the Poisson kernel $P(x, z_{o})).$ The idea now is that there can be only a bounded number $n_{o}$ of points satisfying this property, with $n_{o}$ depending only on the ratio $c_{1}/\varepsilon_{o}.$

 Thus we claim that there is a constant $K_{1} <  \infty ,$ depending only on the choice of $\frac{1}{10}$ in (2.40), and points $p_{1}\in\partial B$ such that
\begin{equation} \label{e2.42}
dist(p_{1},y_{o}) \leq  K_{1}\cdot \rho (y_{o}), 
\end{equation}
and
\begin{equation} \label{e2.43}
m(B_{p_{1}}(\rho (p_{1}))) \leq  \delta_{2}\cdot \varepsilon_{o}. 
\end{equation}
To see this, choose first $p'\in\partial B(1),$ as close as possible to $y_{o}$ such that 
\begin{equation} \label{e2.44}
B_{p'}(\tfrac{1}{2}\rho (p' ))\cap B_{y_{o}}(\rho (y_{o})) = \emptyset  . 
\end{equation}
Note that in general
$$\rho (p' ) \leq  dist(p' , y_{o}) + \rho (y_{o}), $$
so that (2.44) implies
$$dist(p' , y_{o}) \leq  3\rho (y_{o}), $$
and thus
\begin{equation} \label{e2.45}
\rho (p' ) \leq  4\rho (y_{o}). 
\end{equation}
If $p' $ satisfies (2.43), then set $p_{1} = p' .$ If not, so $m(B_{p'}(\rho (p' ))) >  \delta_{2}\cdot \varepsilon_{o},$ then repeat this process with $p' $ in place of $y_{o}.$ Since the total mass is $\delta_{2},$ this can be continued only a bounded number $K_{1}$ of times.

 Clearly,
$$C^{-1}\rho (y_{o}) \leq  \rho (p_{1}) \leq  C\cdot \rho (y_{o}), $$
where $C = C(K_{1}).$ Now choose $y_{1}\in B_{x}(1)\cap B_{p_{1}}(\rho (p_{1})),$ say with $t(y_{1}) = \frac{1}{2}\rho (p_{1}).$ For such a choice, we then have
$$s(y_{1}) \leq  \tfrac{1}{10}\delta_{2}\rho (y_{1})^{-2}, $$
and the result follows.
{\endproof}

 Lemma 2.6 also completes the proof of Theorem 2.5. Finally consider the complementary case to Theorem 2.5. This situation is handled by the following result, which shows that the assumption (2.20) must hold. This result generalizes [An1, Thm.6.2].
\bbgin{theorem} \label{t 2.7.}
  Let $(N, g)$ be a complete ${\cal R}^{2}$ solution with non-negative scalar curvature and uniformly bounded curvature. Then
\begin{equation} \label{e2.46}
liminf_{t\rightarrow\infty} \ \rho^{2}s = 0. 
\end{equation}

\end{theorem}

 The proof of Theorem 2.7 will proceed by contradiction in several steps. Thus, we assume throughout the following that there is some constant $d_{o} > $ 0 such that, for all $x\in (N, g)$,
\begin{equation} \label{e2.47}
s(x) \geq  d_{o}\cdot \rho (x)^{-2}. 
\end{equation}

 Note first that for any complete non-flat manifold, $\rho (x) \leq  2t(x)$ for $t(x)$ sufficiently large, (as in (2.11)), so that (2.47) implies, for some $d > $ 0,
\begin{equation} \label{e2.48}
s(x) \geq  d\cdot  t(x)^{-2}. 
\end{equation}
For reasons to follow later, we assume that $N$ is simply connected, by passing to the universal cover if it is not. Note that (2.47) also holds on any covering space, (with a possibly different constant, c.f. [An1, (3.8)]).

 Consider the conformally equivalent metric
\begin{equation} \label{e2.49}
\Roof{g}{\widetilde} = s\cdot  g. 
\end{equation}
The condition (2.48) guarantees that $(N, \Roof{g}{\widetilde})$ is complete. (Recall that by the maximum principle, $s > $ 0 everywhere, so that $\Roof{g}{\widetilde}$ is well-defined). A standard computation of the scalar curvature $\Roof{s}{\widetilde}$ of $\Roof{g}{\widetilde},$ c.f. [B, Ch.1J] or [An1, (5.18)], gives
\begin{equation} \label{e2.50}
\Roof{s}{\widetilde} = 1 +{\tfrac{2}{3}}\frac{|r|^{2}}{s^{2}} + {\tfrac{3}{2}}\frac{|\nabla s|^{2}}{s^{3}} \geq  1. 
\end{equation}
Thus, $\Roof{g}{\widetilde}$ has uniformly positive scalar curvature.

 We claim that $(N, \Roof{g}{\widetilde})$ has uniformly bounded curvature. To see this, from formulas for the behavior of curvature under conformal changes, c.f. again [B, Ch.1J], one has
\begin{equation} \label{e2.51}
|\Roof{r}{\widetilde}|_{\Roof{g}{\widetilde}} \leq  c_{1}\frac{|r|}{s} + c_{2}\frac{|D^{2}s|}{s^{2}} + c_{3}\frac{|\nabla s|^{2}}{s^{3}}, 
\end{equation}
for some absolute constants $c_{i}.$ Here, the right side of (2.51) is w.r.t. the $g$ metric. The terms on the right on $(N, g)$ are all scale-invariant, so we may estimate them at a point $x\in N$ with $g$ scaled so that $\rho (x, g) = 1$. By assumption (2.47), it follows that $s$ is uniformly bounded below in $B(\frac{1}{2})$ = $B_{x}(\frac{1}{2})$. By Lemma 2.3, $|r|,$ and so also $s$, is uniformly bounded above in $B(\frac{1}{4}).$ Similarly, elliptic regularity for ${\cal R}^{2}$ solutions on $B(\frac{1}{2})$ implies that $|D^{2}s|$ and $|\nabla s|^{2}$ are bounded above on $B(\frac{1}{4}).$ Hence the claim follows.

 Now a result of Gromov-Lawson, [GL, Cor. 10.11] implies, since $N$ is simply connected, that the 1-diameter of $(N, \Roof{g}{\widetilde})$ is at most $12\pi .$ More precisely, let $\Roof{t}{\tilde}(x) = dist_{\Roof{g}{\tilde}}(x, x_{o}).$ Let $\Gamma  = N/\sim ,$ where $x \sim  x' $ if $x$ and $x' $ are in the same arc-component of a level set of $\Roof{t}{\tilde}.$ Then $\Gamma $ is a locally finite metric tree, for which the projection $\pi : N \rightarrow  \Gamma $ is distance non-increasing, c.f. [G, App.1E]. The Gromov-Lawson result states that the diameter of any fiber $F(x) =\pi^{-1}(x)$ in $(N, \Roof{g}{\widetilde})$ is at most $12\pi .$

 Since the curvature of $\Roof{g}{\widetilde}$ is uniformly bounded, it follows that the area of these fibers is also uniformly bounded. In particular, for any $x$,
\begin{equation} \label{e2.52}
vol_{\Roof{g}{\tilde}}B_{x}(r) \leq  C_{o}\cdot  L_{x}(r), 
\end{equation}
where $L(r)$ is the length of the $r$-ball about $x$ in $\Gamma .$

 Further, the uniform curvature bound on $\Roof{g}{\widetilde}$ implies there is a uniform bound $Q$ on the number of edges $E \subset  \Gamma $ emanating from any vertex $v\in\Gamma .$ Note that the distance function $\Roof{t}{\tilde}$ gives $\Gamma $ the structure of a directed tree. Then by construction, at any edge $E \subset\Gamma $ terminating at a vertex $v$, either $\Gamma $ terminates, or there are at least two new edges initiating at $v$. 

 Observe that any point $e$ in an edge $E$, (say not a vertex), divides $\Gamma $ into two components, the inward and outward, with the former containing the base point $x_{o}$ and the latter its complement $\Gamma_{e}.$ The outgoing subtree or branch $\Gamma_{e}$ may have infinite length, in which case it gives an end of $N$, or may have finite length.

 The remainder of the proof needs to be separated into several parts according to the complexity of the graph $\Gamma .$
\bbgin{lemma} \label{l 2.8.}
  Suppose that (2.47) holds. Then the graph $\Gamma $ must have an infinite number of edges. In particular, $N$ has an infinite number of ends. 
\end{lemma}
{\bf Proof:}
 Suppose that $\Gamma $ has only a finite number of edges. It follows that $\Gamma $ has at most linear growth, i.e.
\begin{equation} \label{e2.53}
L_{x}(r) \leq  C_{1}\cdot  r, 
\end{equation}
for some fixed $C_{1} <  \infty .$ Returning to the metric $g$ on $N$, we have $g \leq  C_{2}\cdot  t^{2}\Roof{g}{\widetilde}.$ This together with (2.52) and (2.53) imply that
\begin{equation} \label{e2.54}
vol_{g}(B_{x_{o}}(r)) \leq  C_{3}\cdot  r^{3}, 
\end{equation}
and so Proposition 2.2 implies that $(N, g)$ is flat. This of course contradicts (2.47).

  Since the graph $\Gamma$ is a metric tree, and any vertex has at least two outgoing edges, it follows that $\Gamma$ has infinitely many ends. Hence, by construction, $N$ also has infinitely many ends.

{\endproof}

 The preceding argument will be generalized further in the following to handle the situation when $\Gamma $ has infinitely many edges. To do this however, we first need the following preliminary result. Let $A_{c} = A_{c}(r)$ be any component of the annulus $\Roof{t}{\tilde}^{-1}(r, r+1)$ in $(N, \Roof{g}{\widetilde}).$
\bbgin{lemma} \label{l 2.9.}
   Assume the hypotheses of Theorem 2.7 and (2.47).

{\bf (i).}
 There exists a fixed constant $d <  \infty $ such that
\begin{equation} \label{e2.55}
sup_{A_{c}}s \leq  d\cdot  inf_{A_{c}}s. 
\end{equation}

{\bf (ii).}
 Let $e\in\Gamma ,$ be any edge for which the outgoing branch $\Gamma_{e}$ is infinite, giving an end $N_{e}$ of N. Then
\begin{equation} \label{e2.56}
inf_{N_{e}}s = 0. 
\end{equation}

\end{lemma}
{\bf Proof: (i).}
 We work on the manifold $(N, \Roof{g}{\widetilde})$ as above. From standard formulas for the behavior of $\Delta $ under conformal changes, c.f. [B, Ch.1J], we have
$$\Roof{\Delta}{\widetilde}s = \frac{1}{s}\Delta s + \frac{1}{2s^{2}}|\nabla s|^{2} = \frac{1}{s}\Delta s + \frac{1}{2s}|\Roof{\nabla}{\widetilde}s|^{2}, $$
where the last term is the norm of the gradient, both in the $\Roof{g}{\widetilde}$ metric. Since
$$\Roof{\Delta}{\widetilde}(s^{1/2}) = {\tfrac{1}{2}}s^{-1/2}\bigl(\Roof{\Delta}{\widetilde}s -  \frac{1}{2s}|\Roof{\nabla}{\widetilde}s|^{2}\bigr) , $$
we obtain
$$\Roof{\Delta}{\widetilde}(s^{1/2}) = \tfrac{1}{2}s^{-3/2}\Delta s = -  \tfrac{1}{6}s^{-3/2}|r|^{2} <  0. $$

 Since $(N, \Roof{g}{\widetilde})$ has uniformly bounded geometry, the DeGiorgi-Nash-Moser estimate for supersolutions, c.f. [GT, Thm. 8.18] implies that
$$\bigl(\frac{c}{vol_{\Roof{g}{\tilde}}B}\int_{B}s^{p}d\Roof{V}{\tilde}_{g}\bigr)^{1/p} \leq  inf_{B'  }s, $$
for any concentric $\Roof{g}{\widetilde}$-geodesic balls $B'  \subset  B = B(1)$ and $p <  3/2,$ where $c$ depends only on $dist_{\Roof{g}{\tilde}}(\partial B,\partial B' )$ and $p$. In particular for any component $A_{c}=A_{c}(r)$ of the annulus $\Roof{t}{\tilde}^{-1}(r,r+1)$ in $(N,\Roof{g}{\widetilde}),$ one has
$$\bigl(\frac{1}{vol_{\Roof{g}{\tilde}}A_{c}}\int_{A_{c}}s^{p}dV_{\Roof{g}{\tilde}}\bigr)^{1/p} \leq  c\cdot  inf_{A_{c}'  }s, $$
where $A_{c}' $ is say of half the width of $A_{c}.$ Converting this back to $(N, g)$ gives, after a little calculation,
$$\bigl(\frac{1}{vol_{g}A_{c}}\int_{A_{c}}s^{3-\varepsilon}dV_{g}\bigr)^{1/(3-\varepsilon )} \leq  c\cdot  inf_{A_{c}'  }s, $$
for any $\varepsilon  > $ 0, with $c = c(\varepsilon ).$ By Lemma 2.3, it follows that the $L^{\infty}$ norm of $|r|$ and hence $s$ is bounded on $A_{c}$ by the $L^{2}$ norm average of $s$, which thus gives (2.55).

{\bf (ii).}
 We first note that, in $(N, g)$,
\begin{equation} \label{e2.57}
liminf_{t\rightarrow\infty} \ s = 0.
\end{equation}
This follows easily from the maximum principle at infinity. Thus, let $\{x_{i}\}$ be any minimizing sequence for $s$, so that $s(x_{i}) \rightarrow  inf_{N }s.$ Since $(N, g)$ has bounded curvature, it is clear that $(\Delta s)(x_{i}) \geq  -\varepsilon_{i},$ for some sequence $\varepsilon_{i} \rightarrow $ 0. The trace equation (0.3) then implies that $|r|^{2}(x_{i}) \rightarrow $ 0 which of course implies $s(x_{i}) \rightarrow $ 0.

 Essentially the same argument proves (2.56). Briefly, if (2.56) were not true, one may take any sequence $x_{i}$ going to infinity in $\Gamma_{e}$ and consider the pointed sequence $(N, g, x_{i}).$ Since the curvature is uniformly bounded, a subsequence converges to a complete limit $(N' , g' , x)$, (passing as usual to sufficiently large covers in case of collapse), so that (2.57) holds on the limit. This in turn implies that (2.56) must also hold on $(N, g)$ itself.
{\endproof}

 We are now in position to understand in more detail the structure of $\Gamma .$
\bbgin{lemma} \label{l 2.10.}
  Under the assumptions of Theorem 2.7 and (2.47), there is a uniform upper bound on the distance between nearest vertices in $\Gamma .$ Thus for any vertex $v\in\Gamma ,$ there exists $v'\in\Gamma ,$ with $v'  >  v$ in terms of the direction on $\Gamma ,$ such that
\begin{equation} \label{e2.58}
dist_{\Gamma}(v, v' ) \leq  D, 
\end{equation}
for some fixed $D <  \infty .$
\end{lemma}
{\bf Proof:}
 This is proved by contradiction, so suppose that there is some sequence of edges $E_{i}$ in $\Gamma $ of arbitrarily long length, or an edge $E$ of infinite length. Let $x_{i}$ be the center point of $E_{i}$ in $N$, or a divergent sequence in $E$ in the latter case, and consider the pointed manifolds $(N, \Roof{g}{\widetilde}, x_{i}).$ This sequence has uniformly bounded curvature, and a uniform lower bound on $(volB_{x_{i}}(1), \Roof{g}{\widetilde}),$ since if the sequence volume collapsed somewhere, $(N, \Roof{g}{\widetilde})$ would have regions of arbitrarily long diameter which have the structure of a Seifert fibered space, (c.f. [An1,Thm.2.10]), contradicting by (2.50) [GL,Cor.10.13]. It follows that a subsequence converges to a complete, non-compact manifold $(N' , \Roof{g}{\widetilde}', x)$ with uniformly bounded 1-diameter, uniformly positive scalar curvature, and within bounded Gromov-Hausdorff distance to a line. 

 Consider the same procedure for the pointed sequence $(N, g_{i}, x_{i}),$ where $g_{i} = \rho (x_{i})^{-2}\cdot  g.$ There are now two possibilities for the limiting geometry of $(N, g_i, x_i)$, according to whether liminf \ $\rho (x_{i})/t(x_{i}) =  0$ or $\rho (x_{i})/t(x_{i}) \geq \mu_o $, for some $\mu_o > 0$, as $i \rightarrow \infty$. For clarity, we separate the discussion into these two cases.

{\bf (a).} Suppose that liminf \ $\rho (x_{i})/t(x_{i}) =  0$ and choose a subsequence, also called $\{x_i\}$ such that $\rho (x_{i})/t(x_{i}) \rightarrow  0$. This is again similar to the situation where (2.16) does not hold. Argueing in exactly the same way as in this part of the proof of Lemma 2.4, it follows that one obtains smooth convergence of a further subsequence to a limit $(\bar N, \bar g, \bar x)$. Here, as before, one must pass to suitable covers of larger and larger domains to unwrap a collapsing sequence. The limit $(\bar N, \bar g)$ is a complete, non-flat ${\cal R}^2$ solution. Note that the assumption (2.47) is scale-invariant, invariant under coverings, (c.f. the statement following (2.48)), and invariant under the passage to geometric limits. Hence the estimate (2.47) holds on $(\bar{N}, \bar{g}).$ Now by construction, the graph $\bar{\Gamma}$ associated to $\bar{N}$ is a single line or edge, with no vertices. The argument in Lemma 2.8 above then proves (2.54) holds on $(\bar{N}, \bar{g})$, so that Proposition 2.2 implies that $(\bar{N},\bar{g})$ is flat, contradicting (2.47).

{\bf (b).}
 Suppose that $\rho(x_i)/t(x_i) \geq \mu_o > 0$, for all $i$. In this case, the center points $x_{i}$ remain within uniformly bounded $g_{i}$-distance to the initial vertex $v_{i}$ of $E_{i}.$ The scalar curvature $s_i$ of $g_i$ goes to infinity in a small $g_i$-tubular neighborhood of $v_{i}.$ However, the curvature of $(N, g_i, x_i)$ is uniformly bounded outside the unit ball $(B_{v}(1), g_i)$ and hence in this region one obtains smooth convergence to an incomplete limit manifold $(\bar N', \bar g', \bar x')$, again passing to sufficiently large finite covers to unwrap any collapse. The limit $(\bar N', \bar g')$ is complete away from a compact boundary, (formed by a neigborhood of $\{v\}$). As in Case (a), the limit is a non-flat ${\cal R}^2$ solution satisfying (2.47), and the associated graph $\bar \Gamma'$ consists of a single ray. Hence (2.56) and (2.57) hold on $\bar N'$. 

  Since $\bar \Gamma'$ is a single ray, the annuli $A$ in Lemma 2.9 are all connected and thus (2.56) and (2.57) imply that $\bar s' \rightarrow 0$ uniformly at infinity in $\bar N'$. Hence there is a compact regular level set $L = \{s = s_{o}\}, s_{o} > $ 0, of $s$ for which the gradient $\nabla s|_{L}$ points {\it  out}  of $U$, for $U = \{s \leq  s_{o}\}.$ Observe that $U$ is cocompact in $\bar N'.$

 Now return to the proof of Proposition 2.2, and the trace equation (0.3). Let $\eta $ be a cutoff function as before, but now with $\eta  \equiv $ 1 on a neighborhood of $L$. Integrate as before, but over $U$ in place of $N$ to obtain
$$\int_{U}\eta|r|^{2} = - 3\int_{U}\eta\Delta s = - 3\int_{U}s\Delta\eta  -  3\int_{L}\eta<\nabla s, \nu>  + \ 3\int_{L}s<\nabla\eta , \nu> , $$
where $\nu $ is the outward unit normal. Since $\eta  \equiv $ 1 near $L$, $<\nabla\eta , \nu>  =$ 0. Since $\nabla s|_{L} $ points out of $U$, $<\nabla s, \nu>  > $ 0. Hence,
$$\int_{U}\eta|r|^{2} \leq  - 3\int_{U}s\Delta\eta . $$
From the bound (2.54) obtained as in the proof of Lemma 2.8, the proof of Proposition 2.2 now goes through without any differences and implies that $(\bar{N}, \bar{g})$ is flat in this case also, giving a contradiction.
{\endproof}

 Next, we observe that a similar, but much simpler, argument shows that there is a constant $D <  \infty $ such that for any branch $\Gamma_{e},$ either
\begin{equation} \label{e2.59}
L(\Gamma_{e}) \leq  D, \ \ {\rm or} \ \  L(\Gamma_{e}) = \infty . 
\end{equation}
For suppose there were arbitrarily long but finite branches $\Gamma_{i},$ with length $L(\Gamma_{i}) \rightarrow  \infty $ as $i \rightarrow  \infty ,$ and starting at points $e_{i}.$ Then the same argument proving (2.56) implies that $inf_{N_{i}}s$ cannot occur at or near $e_{i},$ if $i$ is sufficiently large; here $N_{i}$ is the part of $N$ corresponding to $\Gamma_{i}.$ Thus $inf_{N_{i}}s$ occurs in the interior of $N_{i}.$ Since $\Gamma_{i}$ is finite, and thus $\Gamma_{i}$ and $N_{i}$ are compact, this contradicts the minimum principle for the trace equation (0.3).

\medskip

 Given Lemmas 2.9 and 2.10, we now return to the situation following Lemma 2.8 and assume that $(N, g)$ satisfies the assumptions of Theorem 2.7 and (2.47). We claim that the preceding arguments imply the graph $\Gamma $ of $N$ must have exponential growth, in the strong sense that every branch $\Gamma_{e} \subset  \Gamma $ also has exponential growth. 

 Equivalently, we claim that given any point $e$ in an edge $E\subset\Gamma ,$ there is some vertex $v'\in\Gamma_{e},$ (so $v'  >  e$), within fixed distance $D$ to $e$, such that at least two outgoing edges $e_{1}, e_{2}$ from $v' $ have infinite subbranches $\Gamma_{e_{1}}, \Gamma_{e_{2}} \subset  \Gamma_{e}.$ For if this were not the case, then there exist points $e_{i}\in\Gamma $ and branches $\Gamma_{e_{i}} \subset  \Gamma $ and $D_{i} \rightarrow  \infty ,$ such that all subbranches of $\Gamma_{e_{i}}$ starting within distance $D_{i}$ to $e_{i}$ are finite. By (2.59), all subbranches then have a uniform bound $D$ on their length. It follows that one generates as previously in the proof of Lemma 2.10 a geometric limit graph $\Gamma $ with at most linear growth which, as before, gives a contradiction.

 Hence all branches $\Gamma_{e}$ have a uniform rate of exponential growth, i.e. for all $r$ large,
\begin{equation} \label{e2.60}
L(B_{e}(r)) \geq  e^{d\cdot  r}, 
\end{equation}
for some $d > $ 0, where $B_{e}(r)$ is the ball of radius $r$ in $\Gamma_{e}$ about some point in $e$. Since $(N, \Roof{g}{\widetilde})$ has a uniform lower bound on its injectivity radius, $(N_{e}, \Roof{g}{\widetilde})$ satisfies
\begin{equation} \label{e2.61}
vol_{\Roof{g}{\tilde}}(B_{e}(r)) \geq  e^{d\cdot  r}, 
\end{equation}
for some possibly different constant $d$. Further, since $g \geq  c\cdot \Roof{g}{\widetilde},$ (since $s$ is bounded above), (2.61) also holds for $g$, i.e.
\begin{equation} \label{e2.62}
vol_{g}(B_{e}(r)) \geq  e^{d\cdot  r}. 
\end{equation}
again for some possibly different constant $d > 0$.

 However, by (2.57) and (2.55), (and the maximum principle for the trace equation), we may choose a branch $\Gamma_{e}$ such that $s \leq  \delta $ on $N_{e},$ for any prescribed $\delta  > $ 0. Lemma 2.3 then implies that
\begin{equation} \label{e2.63}
|r| \leq  \delta_{1} = \delta_{1}(\delta ), 
\end{equation}
everywhere on $N_{e}.$ From standard volume comparison theory, c.f. [P, Ch.9], (2.63) implies that the volume growth of $N_{e}$ satisfies
\begin{equation} \label{e2.64}
vol_{g}(B_{e}(r)) \leq  e^{\delta_{2}\cdot  r}, 
\end{equation}
where $\delta_{2}$ is small if $\delta_{1}$ is small. Choosing $\delta $ sufficiently small, this contradicts (2.62). This final contradiction proves Theorem 2.7.
{\endproof}

 Theorems 2.5 and 2.7 together prove Theorem 0.1.
\section{Regularity and Apriori Estimates for ${\cal R}_{s}^{2}$ Solutions.}
\setcounter{equation}{0}

 In this section we prove the interior regularity of weak ${\cal R}_{s}^{2}$ solutions as well as apriori estimates for families of ${\cal R}_{s}^{2}$ solutions. These results will be needed in \S 4, and also in [An4]. 

 These are local questions, so we work in a neighborhood of an arbitrary point $x\in N.$ We will assume that $(N, g, x)$ is scaled so that $\rho (x) =$ 1, passing to the universal cover if necessary if the metric is sufficiently collapsed in $B_{x}(\rho (x)).$ In particular, the $L^{2,2}$ geometry of $g$ is uniformly controlled in $B = B_{x}(1).$

 The proof of regularity is similar to the proof of the smooth regularity of ${\cal R}^{2}$ solutions in [An1,\S 4] or that of critical metrics of $I_{\varepsilon},$ for any given $\varepsilon  > $ 0, (c.f. \S 1) in [An1, \S 8]. The proof of apriori estimates proceeds along roughly similar lines, the main difference being that $\alpha $ is not fixed, but can vary over any value in [0, $\infty ).$ Thus one needs uniform estimates, independent of $\alpha .$ To handle this, especially when $\alpha $ is small, one basically uses the interaction of the terms $\alpha\nabla{\cal R}^{2}$ and $L^{*}\omega .$ 

 We assume that $g$ is an $L^{2,2}$ metric satisfying the ${\cal R}_{s}^{2}$ equations
\begin{equation} \label{e3.1}
\alpha\nabla{\cal R}^{2} + L^{*}(\omega ) = 0, 
\end{equation}
\begin{equation} \label{e3.2}
\Delta\omega  = -\frac{\alpha}{4}|r|^{2}, 
\end{equation}
weakly in $B$, with scalar curvature $s \equiv $ 0 (in $L^{2})$ and potential $\omega\in L^{2}.$ If $\alpha  =$ 0, we have assumed, c.f \S 0, that $\omega $ is not identically 0 on $B$. 
\bbgin{theorem} \label{t 3.1.}
  Let (g, $\omega )$ be a weak solution of the ${\cal R}_{s}^{2}$ equations. Then $g$ and $\omega $ are $C^{\infty}$ smooth, in fact real-analytic, in B.
\end{theorem}

 For clarity, the proof will proceed in a sequence of Lemmas. These Lemmas will hold on successively smaller concentric balls $B \supset  B(r_{1}) \supset  B(r_{2}),$ etc, whose ratio $r_{i+1}/r_{i}$ is a definite but arbitrary constant $< $ 1, but close to 1. (The estimates will then depend on this ratio). In other words, we are only considering the interior regularity problem. To simplify notation, we will ignore explicitly stating the size of the ball at each stage and let $\bar{B}$ denote a suitable ball in $B$, with $d = dist(\partial\bar{B}, \partial B).$ Further, $c$ will always denote a constant independent of $\alpha ,$ which is either absolute, or whose dependence is explicitly stated. The value of $c$ may change from line to line or even from one inequality to the next.

 We recall the Sobolev embedding theorem, (in dimension 3), c.f [Ad, Thm. 7.57],
\begin{equation} \label{e3.3}
L^{k,p} \subset  L^{m,q}, \ \ {\rm provided} \ \ k - m <  3/p \ \ {\rm and} \ \ 1/p - (k-m)/3 <  1/q, 
\end{equation}

$$L^{1} \subset  H^{t-2}, \ \ {\rm any} \ \ t <  \frac{1}{2}. $$
where $H^{t},$ (resp. $H_{o}^{t}), t > $ 0, is the Sobolev space of functions (of compact support) with '$t$' derivatives in $L^{2}$ and $H^{-t}$ is the dual of $H_{o}^{t},$ c.f. [Ad, Thm.3.10], [LM, Ch.1.12]. Let $\omega_{o}$ be defined by 
\begin{equation} \label{e3.4}
\omega_{o} = ||\omega||_{L^{2}(B)}. 
\end{equation}

\bbgin{lemma} \label{l 3.2.}
   There is a constant $c =$ c(d) such that
\begin{equation} \label{e3.5}
||\alpha|r|^{2}||_{L^{1}(\bar{B})} \leq  c\cdot \omega_{o}. 
\end{equation}

\end{lemma}
{\bf Proof:}
 This follows immediately from the trace equation (3.2), by pairing it with a suitable smooth cutoff function $\eta $ of compact support in $B$, with $\eta  \equiv $ 1 on $\bar{B}$ and using the self-adjointness of $\Delta .$ Since the $L^{2}$ norm of $\Delta\eta $ is bounded, the left side of (3.2) then becomes bounded by the $L^{2}$ norm of $\omega .$ (This argument is essentially the same as the proof of (2.13)).

{\endproof}
\bbgin{lemma} \label{l 3.3.}
   For any $t <  \frac{1}{2}, \omega\in H^{t}(\bar{B})$, and there is a constant $c =$ c(t, d) such that
\begin{equation} \label{e3.6}
||\omega||_{H^{t}(\bar B)} \leq  c\cdot \omega_{o}. 
\end{equation}

\end{lemma}
{\bf Proof:}
 This also follows from the trace equation (3.2). Namely, Lemma 3.2 implies that the right side of (3.2) is bounded in $L^{1}$ by $\omega_{o}.$ By Sobolev embedding, $L^{1}\subset  H^{t-2},$ for any $t <  \frac{1}{2},$ and we may consider the Laplacian as an operator $\Delta : H^{t} \rightarrow  H^{t-2},$ c.f. [LM, Ch. 2.7] and also [An1, \S 4]. Elliptic theory for $\Delta $ implies that $||\omega||_{H^{t}}$ is bounded by the $H^{t-2}$ norm, and thus $L^{1}$ norm, of the right side and the $L^{2}$ norm of $\omega ,$ each of which is bounded by $\omega_{o}.$

{\endproof}
\bbgin{lemma} \label{l 3.4.}
   The tensor $\alpha r$ is in $H^{t}(\bar{B})$ and there is a constant $c =$ c(t, d) such that
\begin{equation} \label{e3.7}
||\alpha r||_{H^{t}(\bar B)} \leq  c\cdot \omega_{o}. 
\end{equation}

\end{lemma}
{\bf Proof:}
 This follows exactly the arguments of [An1, \S 4], so we will be brief. Write equation (3.1) as
$$D^{*}D\alpha r = Q, $$
where $Q$ consists of all the other terms in (3.1). Using (3.5) and (3.6), it follows that $D^{*}D\alpha r$ is bounded in $H^{t-2}$ by $\omega_{o},$ (c.f. [An1, \S 4] for the details), and the ellipticity of $D^{*}D$ as in Lemma 3.3 above gives the corresponding bound on the $H^{t}$ norm of $\alpha r.$

{\endproof}

From Sobolev embedding (3.3), it follows that
\begin{equation} \label{e3.8}
||\alpha r||_{L^{3-\mu}} \leq  c\omega_{o}, 
\end{equation}
for any given $\mu  = \mu (t) > $ 0. Here and in the following, the norms are on balls $\bar B$, possibly becoming succesively smaller. Using the H\"older inequality, we have
\begin{equation} \label{e3.9}
\int (\alpha|r|^{2})^{p} \leq  \bigl(\int (\alpha|r|)^{pr}\bigr)^{1/r}\bigl(\int|r|^{pq}\bigr)^{1/q}, 
\end{equation}
so that choosing $pr = 3-\mu $ and $pq = 2$ implies that
\begin{equation} \label{e3.10}
||\alpha r^{2}||_{L^{p}} \leq  c\omega_{o}, 
\end{equation}
for any $p = p(t) <  6/5.$ We may now repeat the arguments in Lemma 3.3 above, using the improved estimate (3.10), over the $L^{1}$ estimate. As before, the trace equation now gives, for $p <  6/5,$
$$||\omega||_{L^{1,2-\mu}} \leq  ||\omega||_{L^{2,p}} \leq  c\omega_{o}, $$
for any $\mu  = \mu (p) > $ 0. Now repeat the argument in Lemma 3.4, considering $D^{*}D$ as an operator $L^{1,2-\mu} \rightarrow  L^{-1,2-\mu}.$ This gives
$$||\alpha r||_{L^{6-\varepsilon}} \leq  c||\alpha r||_{L^{1,2-\mu}} \leq  c\omega_{o}, $$
where the first inequality follows from the Sobolev inequality. From the Holder inequality (3.9) again, with $pr = 6-\varepsilon $ and $pq =$ 2, we obtain
$$||\alpha r^{2}||_{L^{p}} \leq  c\omega_{o}, $$
for $p <  3/2.$ Repeating this process again gives
\begin{equation} \label{e3.11}
||\alpha r||_{L^{p}} \leq  c\omega_{o}, 
\end{equation}
for $p <  \infty$ , $c = c(p)$.

 Finally, $\alpha r$ bounded in $L^{p}$ implies $\alpha|r|^{2}$ bounded in $L^{2-\varepsilon}, \varepsilon  = \varepsilon (p),$ so repeating the process once more as above gives $\alpha r$ bounded in $L^{1,6-\varepsilon}$ and bounded in $C^{\gamma}, \gamma  <  \frac{1}{2}.$ This in turn implies $\alpha|r|^{2}$ is bounded in $L^{2}.$ These arguments thus prove the following:
\bbgin{corollary} \label{c 3.5.}
  On a given $\bar{B} \subset  B$, the following estimates hold, with $c = c(d)$:
\begin{equation} \label{e3.12}
||\alpha r||_{L^{1,6}} \leq  c\cdot \omega_{o}, ||\alpha|r|^{2}||_{L^{2}} \leq  c\cdot \omega_{o}, ||\omega||_{L^{2,2}} \leq  c\cdot \omega_{o}. 
\end{equation}
\end{corollary}

{\endproof}

In particular, from Sobolev embedding it follows that $\omega $ is a $C^{1/2}$ function, with modulus of continuity depending on $\omega_{o}.$

 Given these initial estimates, it is now straightforward to prove Theorem 3.1.

\noindent
{\bf Proof of Theorem 3.1.}

 Suppose first $\alpha  > $ 0. Then the iteration above may then be continued indefinitely. Thus since $\alpha r$ is bounded in $L^{1,6} \subset  C^{1/2}$ by Sobolev embedding, $\alpha|r|^{2}$ is also bounded in $L^{1,6}.$ Applying the iteration above, it follows that $\omega\in L^{3,6},$ so that $D^{*}D\alpha r\in L^{1,6}$ implying that $\alpha r\in L^{3,6},$ so that $\alpha|r|^{2}\in L^{3,6},$ and so on. Note that at each stage the regularity of the metric $g$ is improved, since if $\alpha r\in L^{k,p},$ then $g\in L^{k+2,p},$ c.f. [An1, \S 4]. Continuing in this way gives the required $C^{\infty}$ regularity.

 The equations (3.1)-(3.2) form an elliptic system for $(g, \omega)$ with coefficients depending real analytically on $g$. This implies that $g$ and $\omega$ are in fact real-analytic, c.f. [M, Ch.6.6,6.7].

 If $\alpha  =$ 0, then the equations (3.1)-(3.2) are the static vacuum Einstein equations. By assumption the potential $\omega $ is not identically 0 and is an $L^{2,2}$ function by (3.12). The smooth (or real-analytic) regularity of $(g, \omega )$ then follows from standard regularity results for Einstein metrics, c.f. [B, Ch.5].
{\endproof}

 We now turn to higher order estimates for families of ${\cal R}_{s}^{2}$ solutions on balls on which one has initial $L^{2,2}$ control of the metric. This corresponds to obtaining uniform estimates as above which are independent of $\alpha ,$ or equivalently to the smooth compactness of a family of such solutions. In contrast to the ${\cal R}^{2}$ equations, note that the ${\cal R}_{s}^{2}$ equations, for a fixed $\alpha ,$ are not scale-invariant. The full family of ${\cal R}_{s}^{2}$ equations, depending on $\alpha ,$ is scale-invariant. However, $\alpha $ itself is not scale-invariant. As seen in \S 1, $\alpha $ scales inversely to the curvature, i.e. as the square of the distance. On the other hand, the potential $\omega $ is scale-invariant.

 Thus, we now assume we have a sequence or family of smooth ${\cal R}_{s}^{2}$ solutions, with no apriori bound on $\alpha\in [0,\infty )$ or on the $L^{2}$ norm of $\omega $ on $B$. We assume throughout that the metrics are defined on a ball $B = B(1)$, with $r_{h}(B) \geq $ 1. Recall the definition of the constant $c_{o}$ preceding (2.1) in the definition of $r_{h}$ and $\rho ,$ which measures the deviation of $g$ from the flat metric. We define the $L^{k,2}$ curvature radius $\rho^{k}$ exactly in the same way as the $L^{2}$ curvature radius, with $|\nabla^{k} r|$ in place of $|r|$, c.f. [An1, Def.3.2]. 
\bbgin{theorem} \label{t 3.6}
 Let (N, g, $\omega )$ be an ${\cal R}_{s}^{2}$ solution, with $x\in N$ and $r_{h}(x) \geq $ 1. Suppose
\begin{equation} \label{e3.13}
\omega  \leq  0, 
\end{equation}
in $B_{x}(1).$ Then there is a constant $\rho_{o} > $ 0 such that the $L^{k,2}$ curvature radius $\rho^{k}$ satisfies, for all $k\geq 1,$
\begin{equation} \label{e3.14}
\rho^{k}(x) \geq  \rho_{o}. 
\end{equation}
In particular, by Sobolev embedding, the Ricci curvature $r$ is bounded in $C^{k}.$ Further the $L^{k+2,2}$ norm of the potential $\omega $ is uniformly bounded in $B_{x}(\rho_{o}),$ in that, for a given $\varepsilon  > $ 0, which may be made (arbitrarily) small if $c_{o}$ is chosen sufficiently small,
\begin{equation} \label{e3.15}
||\omega||_{L^{k+2,2}(B_{x}(\rho_{o}))} \leq  max\bigl( c(\varepsilon^{-1},k)\cdot ||\omega||_{L^{1}(B_{x}(\varepsilon ))}, c_{o}\bigr) , 
\end{equation}

\end{theorem}

 The proof of this result will again be carried out in a sequence of steps. This result is considerably more difficult to prove than Theorem 3.1, since one needs to obtain estimates independent of $\alpha $ and $\omega_{o}.$ To do this, we will strongly use the assumption (3.13) that $\omega  \leq $ 0. It is not clear if Theorem 3.6 holds in general when $\omega  \geq $ 0. In fact, this is the main reason why the hypothesis (3.13) is needed in Theorem 0.2.

 We note that by Theorem 3.1, the regularity of $(g, \omega )$ is not an issue, so we will assume that $g$ and $\omega $ are smooth, (or real-analytic) on $B$. 

 Theorem 3.6 is easy to prove from the preceding estimates if one has suitable bounds on $\alpha $ and $\omega_{o}.$ It is worthwhile to list these explicitly.

{\bf (i).}
 Suppose $\alpha $ is uniformly bounded away from 0 and $\infty ,$ and $\omega_{o}$ is bounded above, i.e.
\begin{equation} \label{e3.16}
\kappa  \leq  \alpha  \leq  \kappa^{-1}, \omega_{o} \leq  \kappa^{-1} 
\end{equation}
for some fixed $\kappa  > $ 0. One may then just carry out the arguments preceding Corollary 3.5 and in the proof of Theorem 3.1 to prove Theorem 3.6, with bounds then depending on $\kappa .$ This gives (3.14) and also (3.15), but with $\omega_{o}$ in place of the $L^{1}$ norm of $\omega $ in $B_{x}(\varepsilon ).$ In case 
$$||\omega||_{L^{1}(B_{x}(\varepsilon ))}<<  \omega_{o}, $$
we refer to Case (I) below, where this situation is handled.

{\bf (ii).}
 Suppose that $\alpha $ is uniformly bounded away from 0, and $\omega_{o}$ is bounded above, i.e.
\begin{equation} \label{e3.17}
\kappa  \leq  \alpha , \omega_{o} \leq  \kappa^{-1} 
\end{equation}
for some fixed $\kappa  > $ 0. Thus, the difference with Case (i) is that $\alpha $ may be arbitrarily large. Then the equations (3.1)-(3.2) may be renormalized by dividing by $\alpha ,$ so that $\alpha $ becomes 1. This decreases the $L^{2}$ norm of the potential $\omega ,$ i.e. $\omega_{o},$ but otherwise leaves the preceding arguments unchanged. Hence, Theorem 3.6 is again proved in this case, with bounds depending on $\kappa .$ Observe that both of these arguments in (i), (ii) do not require the bound (3.13).

\medskip

 The main difficulty is when $\alpha $ is small, and especially when $\omega_{o}$ is small as well. In this situation, the equations (3.1)-(3.2) approach merely the statement that 0 $\sim $ 0. Here one must understand the relative sizes of $\alpha $ and $\omega_{o}$ to proceed further.

 First note that by Theorem 3.1, we may assume that $\omega$ is not identically zero on any ball $\bar{B} \subset B$, since otherwise, (by real-analyticity), $\omega \equiv 0$ on $B$, and hence the solution is flat; (recall that we have assumed $\alpha > 0$ if $\omega \equiv 0$ on $B$). It follows that the set $\{\omega  =$ 0\} is a closed real-analytic set in B.

 To begin, we must obtain apriori estimates on, for example, the $L^{2}$ norm of $\omega $ in terms of its value at or near the center point $x$, c.f. Corollary 3.9. The first step in this is the following $L^{1}$ estimate, which may also be of independent interest. (Of course this is meant to apply to $u = -\omega $ and $f = \frac{1}{4}\alpha|r|^{2}$, as in (3.2)). 
\bbgin{proposition} \label{p 3.7.}
   Let (B, g), $B = B_{x}(1)$ be a geodesic ball of radius 1, with $r_{h}(x) \geq $ 1. For a given smooth function $f$ on B, let $u$ be a solution of
\begin{equation} \label{e3.18}
\Delta u = f, 
\end{equation}
on B, with
\begin{equation} \label{e3.19}
u \geq  0.
\end{equation}
Let $B_{\varepsilon} = B_{x}(\varepsilon ),$ where $\varepsilon  > $ 0 is (arbitrarily) small, depending only on the choice of $c_{o}.$ Then there is a constant $c > $ 0, depending only on $\varepsilon^{-1},$ such that
\begin{equation} \label{e3.20}
\int_{B}u \leq  c(\varepsilon^{-1})\int_{B}|f| + {\tfrac{1}{2}} u_{av}(\varepsilon ) + \int_{B_{\varepsilon}}u, 
\end{equation}
where $u_{av}$ is the average of $u$ on $S_{\varepsilon} = \partial B_{\varepsilon}.$
\end{proposition}
{\bf Proof:}
 Let $\eta $ be a non-negative function on $B$ such that $\eta  = |\nabla\eta| =$ 0 on $\partial B,$ determined more precisely below. Multiply (3.18) by $\eta $ and apply the divergence theorem to obtain
\begin{equation} \label{e3.21}
\int_{B \setminus {B_{\varepsilon}}}\eta f = \int_{B \setminus {B_{\varepsilon}}}\eta\Delta u = \int_{B \setminus {B_{\varepsilon}}}u\Delta\eta  + \int_{S_{\varepsilon}}\eta<\nabla u, \nu>  -  \int_{S_{\varepsilon}}u<\nabla\eta , \nu>,
\end{equation}
where $\nu $ is the outward unit normal. Let $\{x_{i}\}$ be a harmonic coordinate chart on $B$, so that the metric $g$ is bounded in $L^{2,2}$ on $B$ in the coordinates $\{x_{i}\},$ since $r_{h}(x) \geq $ 1. We may assume w.l.o.g. that in these coordinates, $g_{ij}(x) = \delta_{ij}.$ Let $\sigma  = (\sum x_{i}^{2})^{1/2}.$ Then the ratio $\sigma /t,$ for $t(y) = dist_{g}(y,x),$ satisfies
\begin{equation} \label{e3.22}
1- c_{1} \leq  \frac{\sigma}{t} \leq  1+c_{1}, 
\end{equation}
in $B$, where the constant $c_{1}$ may be made small by choosing the constant $c_{o}$ in the definition of $r_{h}$ sufficiently small.

 We choose $\eta  = \eta (\sigma )$ so that $\eta (1) =$ 0 and $\eta' (1) =$ 0 and so that $\Delta\eta $ is close to the constant function 2 in $C^{0}(B \setminus {B_{\varepsilon}}),$ i.e.
\begin{equation} \label{e3.23}
\Delta\eta  \sim  2. 
\end{equation}
(Actually, since $\{\sigma  =$ 1\} may not be contained in $B$, one should replace by the boundary condition at 1 by the same condition at $\sigma  = 1-\delta , \delta  = \delta (c_{o})$ small, but we will ignore this minor adjustment below).

 To determine $\eta$, we have $\Delta\eta  = \eta'\Delta\sigma  + \eta''|\nabla\sigma|^{2}.$ Since the metric $g = g_{ij}$ is close to the flat metric $\delta_{ij},$ we have $||\nabla\sigma|^{2}- 1| <  \delta $ and $|\Delta\sigma-\frac{2}{\sigma}| <  \frac{\delta}{\sigma}$, in $C^{0}$, where $\delta $ may be made (arbitrarily) small by choosing $c_{o}$ sufficiently small. Thus, let $\eta $ be a solution to
\begin{equation} \label{e3.24}
\frac{2}{\sigma}\eta'  + \eta''  = 2. 
\end{equation}
It is easily verified that the function
\begin{equation} \label{e3.25}
\eta (\sigma ) = \tfrac{2}{3}\sigma^{-1} -  1 + \tfrac{1}{3}\sigma , 
\end{equation}
satisfies (3.24), with the correct boundary conditions at $\sigma  =$ 1. Observe that $\eta (\varepsilon ) \sim  \frac{2}{3}\varepsilon^{-1}, \ \eta' (\varepsilon ) \sim  -<\nabla\eta , \nu> (\varepsilon ) \sim  -\frac{2}{3}\varepsilon^{-2},$ for $\varepsilon $ small. Hence, (3.23) is satisfied on $B \setminus {B_{\varepsilon}}$ for a given choice of $\varepsilon $ small, if $\delta ,$ i.e. $c_{o},$ is chosen sufficiently small.

 With this choice of $\eta $ and $\varepsilon ,$ (3.21) becomes
$$\int_{B \setminus {B_{\varepsilon}}}\eta f \sim  2\int_{B \setminus {B_{\varepsilon}}}u  -  {\tfrac{2}{3}} \varepsilon^{-2}\int_{S_{\varepsilon}}u  + {\tfrac{2}{3}} \varepsilon^{-1}\int_{S_{\varepsilon}}<\nabla u, \nu> . $$
Applying the divergence theorem on $B_{\varepsilon}$ to the last term gives
$$ {\tfrac{2}{3}} \varepsilon^{-1}\int_{S_{\varepsilon}}<\nabla u, \nu>  = -{\tfrac{2}{3}} \varepsilon^{-1}\int_{B_{\varepsilon}}\Delta u  \sim  - c\cdot \varepsilon^{-1}\int_{B_{\varepsilon}}f. $$
Thus, combining these estimates gives
$$ 2\int_{B \setminus {B_{\varepsilon}}}u \leq   c(\varepsilon^{-1})\int_{B}|f| + u_{av}(\varepsilon ), $$
which implies (3.20).

{\endproof}

 The next Lemma is a straightforward consequence of this result.
\bbgin{lemma} \label{l 3.8.}
  Assume the same hypotheses as Proposition 3.7. Then for any $\mu  > $ 0, there is a constant $c = c(\mu ,d, \varepsilon^{-1})$ such that
\begin{equation} \label{e3.26}
||u||_{L^{3-\mu}(\bar{B})} \leq  c(\mu ,d,\varepsilon^{-1})\cdot \int_{B}|f| + {\tfrac{1}{2}} u_{av}(\varepsilon ). 
\end{equation}

\end{lemma}
{\bf Proof:}
 Let $\phi $ be a smooth positive cutoff function, supported in $B$, with $\phi  \equiv $ 1 on $\bar{B}.$ As in Lemma 3.3, the Laplacian $\Delta $ is a bounded surjection $\Delta : H^{t} \rightarrow  H^{t-2},$ for $t <  \frac{1}{2},$ so that
\begin{equation} \label{e3.27}
||\phi u||_{H^{t}} \leq  c\cdot ||\Delta (\phi u)||_{H^{t-2}}. 
\end{equation}
By Sobolev embedding  $||\phi u||_{L^{3-\mu}} \leq  c\cdot ||\phi u||_{H^{t}}, c = c(\mu ), \mu  = \mu (t)$ so that it suffices to estimate the right side of (3.27). By definition,
\begin{equation} \label{e3.28}
||\Delta (\phi u)||_{H^{t-2}} = sup|\int h\Delta (\phi u)|, 
\end{equation}
where the supremum is over compactly supported functions $h$ with $||h||_{H_{o}^{2-t}} \leq $ 1. Since $\Delta u = f$, we have 
$$\Delta (\phi u) = \phi f + u\Delta\phi  + 2<\nabla\phi ,\nabla u> . $$
Now
\begin{equation} \label{e3.29}
\int h\phi f \leq  c\cdot ||f||_{L^{1}}\cdot ||h||_{L^{\infty}} \leq  c\cdot ||f||_{L^{1}}, 
\end{equation}
where the last inequality follows from Sobolev embedding. Next
\begin{equation} \label{e3.30}
\int hu\Delta\phi  \leq  c\cdot ||h||_{L^{\infty}}||u||_{L^{1}} \leq  c\cdot ||u||_{L^{1}}, 
\end{equation}
since we may choose $\phi $ with $|\Delta\phi| \leq  c$ (e.g. $\phi  = \phi (\sigma ),$ for $\sigma $ as in the proof of Prop. 3.7). For the last term, applying the divergence theorem gives
\begin{equation} \label{e3.31}
\int h<\nabla\phi , \nabla u>  = -\int uh\Delta\phi  -\int u<\nabla h, \nabla\phi> . 
\end{equation}
The first term here is treated as in (3.30). For the second term, for any $\varepsilon_{1} > $ 0 we may choose $\phi $ so that $|\nabla\phi^{\varepsilon_{1}}|$ is bounded, (with bound depending only on $\varepsilon_{1}).$ By Sobolev embedding $|\nabla h|$ is bounded in $L^{3+\varepsilon_{2}},$ for $\varepsilon_{2} = \varepsilon_{2}(\mu )> $ 0. Thus, applying the H\"older inequality,
$$|\int u<\nabla h, \nabla\phi>| \leq  c\int\phi^{1-\varepsilon_{1}}u|\nabla h| \leq  c\cdot ||\phi^{1-\varepsilon_{1}}u||_{L^{(3/2)-\varepsilon_{2}}}. $$
Write $(\phi^{1-\varepsilon_{1}}u)^{(3/2)-\varepsilon_{2}} = (\phi u)^{p}u^{1/2},$ where $p = 1-\varepsilon_{2}.$ For $\varepsilon_{2}$ small, this gives $\varepsilon_{1} \sim  1/3.$ It follows that 
$$\int (\phi^{1-\varepsilon_{1}}u)^{(3/2)-\varepsilon_{2}} = \int (\phi u)^{p}\cdot  u^{1/2} \leq  \bigl(\int (\phi u)^{2p}\bigl)^{1/2}\bigl(\int u\bigr)^{1/2}. $$
This estimate, together with standard use of the Young and interpolation inequalities, c.f. [GT, (7.5),(7.10)] implies that
\begin{equation} \label{e3.32}
||\phi^{1-\varepsilon_{1}}u||_{L^{(3/2)-\varepsilon_{2}}} \leq  \kappa||\phi u||_{L^{2}} + c\cdot \kappa^{-a}||u||_{L^{1}}, 
\end{equation}
for any given $\kappa  > $ 0, where $c$ and $a$ depend only on $\varepsilon_{2}.$ Since $||\phi u||_{L^{2}} \leq  ||\phi u||_{H^{t}},$ by choosing $\kappa $ small in (3.32), we may absorb the first term on the right in (3.32) to the left in (3.27). Combining the estimates above gives
$$||u||_{L^{3-\mu}(\bar{B})} \leq  c(\mu )\cdot  (||f||_{L^{1}(B)} + ||u||_{L^{1}(B)}). $$
The $L^{1}$ norm of $u$ is estimated in Proposition 3.7. In addition, we have
$$\int_{B_{\varepsilon}}u \leq  ||u||_{L^{3-\mu}(\bar{B})}\cdot  (volB_{\varepsilon})^{q}, $$
where $q = \frac{2-\mu}{3-\mu}.$ Choosing $\varepsilon $ small, this term may also be absorbed into the right in (3.27), which gives the bound (3.26).

{\endproof}

 We point out that although the assumption (3.19) on $u$, (equivalent to (3.13) on $\omega$), is in fact not used in Proposition 3.7, it is required in Lemma 3.8 to control the $L^1$ norm of $u$. 
 We now apply these results to the trace equation (3.2), so that $u = -\omega $ and $f = \frac{1}{4}\alpha|r|^{2}.$ Using the fact that $f \geq $ 0, we obtain:
\bbgin{corollary} \label{c 3.9.}
  On (B, g) as above,
\begin{equation} \label{e3.33}
sup_{\bar{B}}|\omega| \leq  c(d,\varepsilon^{-1})\cdot \int_{B}\alpha|r|^{2} + {\tfrac{1}{2}} |\omega|_{av}(\varepsilon ). 
\end{equation}

\end{corollary}
{\bf Proof:}
 By Lemma 3.8, it suffices to prove there is a bound
$$sup_{\bar{B}}|\omega| \leq  c||\omega||_{L^{2}(B)}. $$
Since $f$ is positive, i.e. $\Delta (-\omega ) \geq $ 0, this estimate is an application of the DeGiorgi-Nash-Moser sup estimate for subsolutions of the Laplacian, c.f. [GT, Thm 8.17].

{\endproof}

 We are now in a position to deal with the relative sizes of $\alpha $ and $\omega .$ First, we renormalize $\omega $ to unit size near the center point $x$. Observe that for any $\varepsilon  > $ 0,
\begin{equation} \label{e3.34}
|\omega|_{av}(\varepsilon ) >  0, 
\end{equation}
since if $|\omega|_{av}(\varepsilon ) =$ 0, then by the minimum principle applied to the trace equation (3.2), $\omega  \equiv $ 0 in $B_{\varepsilon}$ which has been ruled out above. For the remainder of the proof, we fix an $\varepsilon  > $ 0 small and set
\begin{equation} \label{e3.35}
\bar{\omega} = \frac{\omega}{|\omega|_{av}(\varepsilon )}. 
\end{equation}
Similarly, divide the equations (3.1)-(3.2) by $|\omega|_{av}(\varepsilon );$ this has the effect of changing $\alpha $ to
\begin{equation} \label{e3.36}
\bar{\alpha} = \frac{\alpha}{|\omega|_{av}(\varepsilon )}. 
\end{equation}
By Corollary 3.9, (or Lemma 3.8), the $L^{2}$ norm of $\bar{\omega}$ on $\bar{B},$ i.e. $\bar{\omega}_{o},$ is now controlled by the $L^{1}$ norm of $\bar{\alpha}|r|^{2}$ on $B$,
\begin{equation} \label{e3.37}
\bar{\omega}_{o} \leq  c\cdot  [||\bar{\alpha}|r|^{2}||_{L^{1}(B)}+1]. 
\end{equation}

 We now divide the discussion into three cases, similar to the discussion following Theorem 3.6.

{\bf (I).}
 Suppose $\bar{\alpha}$ is uniformly bounded away from 0 and $\infty ,$ i.e.
\begin{equation} \label{e3.38}
\kappa  \leq  \bar{\alpha} \leq  \kappa^{-1}, 
\end{equation}
for some $\kappa  > $ 0. Since $r_{h}(x) \geq $ 1, the $L^{1}$ norm of $\bar{\alpha}|r|^{2}$ on $B$ is then uniformly bounded, and so $\bar{\omega}_{o}$ is uniformly bounded by (3.37). Hence, the proof of (3.14) follows as in Case (i) above, starting with Corollary 3.5. Further, 
$$\bar{\omega}_{o} \leq  c(\varepsilon^{-1})\cdot ||\bar{\omega}||_{L^{1}(B_{\varepsilon})}, $$
since both terms are uniformly bounded away from 0 and $\infty .$ Hence (3.15) also follows as in Case (i) above.

{\bf (II).}
 Suppose $\bar{\alpha}$ is uniformly bounded away from 0, i.e.
\begin{equation} \label{e3.39}
\kappa  \leq  \bar{\alpha}, 
\end{equation}
for some $\kappa  > $ 0. Then as in Case (ii), divide $\bar{\omega}$ and $\bar{\alpha}$ by $\bar{\alpha},$ so that the coefficient $\alpha'  = \frac{\bar{\alpha}}{\bar{\alpha}}  =$ 1. Then $\omega'  = \frac{\bar{\omega}}{\bar{\alpha}}$ has bounded $L^{2}$ norm on $\bar{B}$ by (3.37) and the proof of (3.14) proceeds as in Case (I). In this case, the $L^{1}$ norm of $\omega' $ on $B_{\varepsilon}$ may be very small, but (3.15) follows as above from the renormalization of (3.37).

{\bf (III).}
 $\bar{\alpha}$ is small, i.e.
\begin{equation} \label{e3.40}
\bar{\alpha} \leq  \kappa . 
\end{equation}
Again from (3.37), it follows that $\bar{\omega}_{o}$ is uniformly bounded on $\bar{B}.$ Thus, by Corollary 3.5, we have
\begin{equation} \label{e3.41}
||\bar{\omega}||_{L^{2,2}} \leq  c(\kappa ), 
\end{equation}
so that, by Sobolev embedding, $\bar{\omega}$ is uniformly controlled in the $C^{1/2}$ topology inside $\bar{B}.$ This, together with (3.35) implies that $\bar{\omega}(x) \sim  - 1,$ and hence there is a smaller ball $B'\subset\bar{B},$ whose radius is uniformly bounded below, such that for all $y\in B' ,$
\begin{equation} \label{e3.42}
\bar{\omega}(y) \leq  -\tfrac{1}{10}. 
\end{equation}
(The factor $10^{-1}$ can be replaced by any other fixed negative constant).

\medskip

 We now proceed with the proof of Theorem 3.6 under these assumptions. For notational simplicity, we remove the bar and prime from now on, so that $B' $ becomes $B$, $\bar{\omega}$ becomes $\omega , \bar{\alpha}$ becomes $\alpha ,$ and $\alpha , \omega $ satisfy (3.40)-(3.42) on $B$. The previous arguments, when $\alpha $ is bounded away from 0, essentially depend on the (dominant) $\alpha\nabla{\cal R}^{2}$ term in (3.1), and much less on the $L^{*}(\omega )$ term, (except in the estimates on $\omega $ in Proposition 3.7 - Corollary 3.9). When $\alpha $ is or approaches 0, one needs to work mainly with the $L^{*}(\omega )$ term to obtain higher order estimates, ensuring that the two terms $\alpha\nabla{\cal R}^{2}$ and $L^{*}(\omega )$ do not interfere with each other. For the following Lemmas, the preceding assumptions are assumed to hold.
\bbgin{lemma} \label{l 3.10.}
  There is a constant $c <  \infty $ such that, on $\bar{B} \subset B,$
\begin{equation} \label{e3.43}
\alpha (\int|r|^{6})^{1/3} \leq  c, \ \  \alpha\int|r|^{3} \leq  c, \ \ \alpha\int|Dr|^{2} \leq  c.  
\end{equation}

\end{lemma}
{\bf Proof:}
 For a given smooth cutoff function $\eta $ of compact support in $\bar{B},$ pair (3.1) with $\eta^{2}r$ and integrate by parts to obtain
\begin{equation} \label{e3.44}
\alpha\int|D\eta r|^{2} -  \int\omega\eta^{2}|r|^{2} \leq  \alpha c\int\eta^{2}|r|^{3} + \int< D^{2}(-\omega ), \eta^{2}r>  + \alpha c\int (|D\eta|^{2}+ |\Delta\eta^{2}|)|r|^{2}. 
\end{equation}
By the Sobolev inequality,
\begin{equation} \label{e3.45}
(\int (\eta|r|)^{6})^{1/3} \leq  c_{s}\int|D(\eta|r|)|^{2} \leq  cc_{s}\int||D\eta r||^{2}. 
\end{equation}
where $c_{s}$ is the Sobolev constant of the embedding $L_{0}{}^{1,2} \subset  L^{6}.$ Also, by the H\"older inequality,
\begin{equation} \label{e3.46}
\int\eta^{2}|r|^{3} \leq  (\int (\eta|r|)^{6})^{1/3}\cdot  (\int|r|^{3/2})^{2/3} \leq  c_{o}(\int (\eta|r|)^{6})^{1/3}, 
\end{equation}
where the last inequality follows from the definition of $\rho ,$ c.f. the beginning of \S 2. We may assume $c_{o}$ is chosen sufficiently small so that $cc_{o}c_{s} \leq  \frac{1}{2}.$ Thus, the cubic term on the right in (3.44) may be absorbed into the left, giving
$$\frac{c}{2}\alpha (\int (\eta r)^{6})^{1/3} -  \int\omega\eta^{2}|r|^{2} \leq  \int< D^{2}(-\omega ), \eta^{2}r>  + \alpha c\int (|D\eta|^{2}+|\Delta\eta^{2}|)|r|^{2}. $$
For the first term on the right, we may estimate
$$\int< D^{2}(-\omega ), \eta^{2}r> \  \leq  (\int|D^{2}\omega|^{2})^{1/2}(\int|r|^{2})^{1/2} \leq  c, $$
by (3.41) and the definition of $\rho .$ For the second term, $|D\eta|$ and $|\Delta\eta^{2}|$ are uniformly bounded in $L^{\infty}(B),$ so that
$$\alpha\int (|D\eta|^{2}+ |\Delta\eta^{2}|)|r|^{2} \leq  c, $$
either by Lemma 3.2 or by (3.40) and the definition of $\rho .$ From (3.42), $-\omega  \geq  \frac{1}{10}$ in $\bar{B},$ so that one obtains the bound
$$\alpha (\int|r|)^{6})^{1/3} \leq  c,  $$
and thus also
$$\alpha\int|Dr|^{2} \leq  c. $$
The bound on the $L^{1}$ norm of $\alpha|r|^{3}$ then follows from (3.46).

{\endproof}

\noindent
\bbgin{lemma} \label{l 3.11.}
  There is a constant $c <  \infty $ such that, on $\bar{B} \subset B,$
\begin{equation} \label{e3.47}
\int|r|^{3} \leq  c, \ \  \alpha (\int|r|^{9})^{1/3} \leq  c. 
\end{equation}

\end{lemma}
{\bf Proof:}
 By pairing (3.1) with $r$, one obtains the estimate
\begin{equation} \label{e3.48}
-\alpha\Delta|r|^{2} + \alpha|Dr|^{2}  -  \omega|r|^{2} \  \leq  \ < D^{2}(-\omega ), r>  + \ \alpha|r|^{3}. 
\end{equation}
Multiply (3.48) by $\eta|r|$ and integrate by parts, using the bound (3.42) and the Sobolev inequality as in (3.45) to obtain
\begin{equation} \label{e3.49}
c\cdot \alpha (\int|r|^{9})^{1/3}+ \frac{1}{10}\int|r|^{3} \leq  \alpha\int|r|^{4} + \int< D^{2}(-\omega ), \eta|r|r> + \alpha\int|r|^{2}|D|r||\cdot |D\eta|. 
\end{equation}
We have
$$\int< D^{2}(-\omega ), \eta|r|r> \  \leq  \int|r|^{2}|D^{2}\omega| \leq  (\int|r|^{3})^{2/3}(\int|D^{2}\omega|^{3})^{1/3} ,$$
and
$$(\int|D^{2}\omega|^{3})^{1/3} \leq  c(\int\alpha^{3}|r|^{6})^{1/3}+ c(\int|\omega|^{2})^{1/2} \leq  c, $$
where the last estimate follows from elliptic regularity on the trace equation (3.2) together with Lemma 3.10 and (3.41). Thus, this term may be absorbed into the left of (3.49), (because of the $2/3$ power), (unless $||r||_{L^{3}}$ is small, in which case one obtains a bound on $||r||_{L^{3}}).$ 

 Similarly
$$\alpha\int|r|^{4} = \alpha\int|r|^{3}|r| \leq  \alpha (\int|r|^{9})^{1/3}(\int|r|^{3/2})^{2/3} \leq  c_{o}\alpha (\int|r|^{9})^{1/3}. $$
Since $c_{o}$ is sufficiently small, this term may be absorbed into the left. For the last term on the right in (3.49),
$$\alpha\int|r|^{2}|D|r||\cdot |D\eta| \leq  \delta\alpha\int|r|^{4} + \delta^{-1}\alpha\int|D|r||^{2}|D\eta|^{2};   $$
the second term on the right is bounded by Lemma 3.10, while as in the previous estimate, the first term may be absorbed into the left, by choosing $\delta $ sufficiently small. These estimates combine to give the result.

{\endproof}

 This same argument may be repeated once more, pairing with $\eta|r|^{2},$ to obtain a bound
\begin{equation} \label{e3.50}
\int|r|^{4} \leq  c. 
\end{equation}

\noindent
\bbgin{lemma} \label{l 3.12.}
  There is a constant $c <  \infty $ such that, on $\bar{B} \subset B,$
\begin{equation} \label{e3.51}
\alpha^{1/2}||D^{*}Dr||_{L^{2}} \leq  c,\ \ ||Dr||_{L^{2}} \leq  c, \ \  ||D^{2}\omega||_{L^{1,2}} \leq  c, 
\end{equation}

\end{lemma}
{\bf Proof:}
 We may take the covariant derivative of (3.1) to obtain
\begin{equation} \label{e3.52}
\alpha DD^{*}Dr +  \alpha Dr^{2} -  D\omega r = DD^{2}(-\omega ), 
\end{equation}
where $Dr^{2}$ denotes derivatives of terms quadratic in curvature. Hence
$$\alpha< DD^{*}Dr,\eta Dr>  + \alpha< Dr^{2},\eta Dr>  -  < D\omega r,\eta Dr> \ = -< DD^{2}(\omega ),\eta Dr> . $$
Using (3.42), this gives rise to the bound
\begin{equation} \label{e3.53}
\alpha\int|D^{*}D\eta r|^{2}+\frac{1}{10}\int\eta|Dr|^{2}\leq 
\end{equation}
$$\leq  \alpha c\int\eta|r||Dr|^{2}+\delta\int\eta|Dr|^{2}+c\delta^{-1}\int\eta|DD^{2}\omega|^{2}+\alpha c\int|D\eta|^{2}|Dr|^{2}+c\int\eta|D\omega||Dr||r|. $$
From elliptic regularity together with the Sobolev inequality (3.45), we have, (ignoring some constants),
\begin{equation} \label{e3.54}
\alpha (\int|D\eta r|^{6})^{1/3}\leq  \alpha\int|D^{*}D\eta r|^{2}. 
\end{equation}
Using same argument kind of argument as in (3.46), the first term on right of (3.53) can then be absorbed into left. Similarly, the second term in (3.53) can be absorbed into the $\frac{1}{10}$ term on the left, for $\delta $ small. For the third term, elliptic regularity gives, (ignoring constants and lower order terms),
\begin{equation} \label{e3.55}
\int|DD^{2}\omega|^{2} \leq  \int|D\Delta\omega|^{2}+c \leq  \alpha^{2}\int|r|^{2}|D|r||^{2}+c \leq  \alpha^{2}(\int|r|^{4})^{1/2}(\int|D|r||^{4})^{1/2}+c. 
\end{equation}
But by Lemma 3.10,
$$\alpha (\int|r|^{4})^{1/2}\leq  c\cdot  c_{o}, $$
so that this third term may also be absorbed in left using (3.54). The fourth term on the right in (3.53) is bounded by Lemma 3.10. Finally for the last term, the H\"older inequality gives
$$\int\eta|D\omega||Dr||r| \leq  ||D\omega||_{L^{6}}||r||_{L^{3}}||Dr||_{L^{2}}. $$
The first two terms on the right here are bounded, by (3.41) and Lemma 3.11, so by use of the Young inequality, this term may also be absorbed into the left in (3.53).

 It follows that the left side of (3.53) is uniformly bounded. This gives the first two bounds in (3.51), while the last bound follows from the argument (3.55) above.

{\endproof}

 The estimates in Lemma 3.12, together with the previous work in Cases (I)-(III) above, proves Theorem 3.6 for $k =$ 1. By taking higher covariant derivatives of (3.1) and continuing in the same way for Case(III), one derives higher order estimates on $g$ and $\omega .$ 

This completes the proof of Theorem 3.6.

{\endproof}

\section{Non-Existence of ${\cal R}_{s}^{2}$ Solutions with Free $S^{1}$ Action.}
\setcounter{equation}{0}

 In this section, we prove Theorem 0.2. Thus let $N$ be an open oriented 3-manifold and $g$ a complete Riemannian metric on $N$ satisfying the ${\cal R}_{s}^{2}$ equations (0.4)-(0.5). The triple $(N, g, \omega )$ is assumed to admit a free isometric $S^{1}$ action. The assumption (0.6) on $\omega $ will enter only later.

 Let $V$ denote the quotient space $V = N/S^{1}.$ By passing to a suitable covering space if necessary, we may and will assume that $V$ is simply connected, and so $V = {\Bbb R}^{2}$ topologically.

 The metric $g$ is a Riemannian submersion to a complete metric $g_{V}$ on $V$. Let $f: V \rightarrow  {\Bbb R} $ denote the length of the orbits, i.e. $f(x)$ is the length of the $S^{1}$ fiber through $x$. Standard submersion formulas, c.f. [B, 9.37], imply that that the scalar curvature $s_{V}$ of $g_{V},$ equal to twice the Gauss curvature, is given by
\begin{equation} \label{e4.1}
s_{V} = s + |A|^{2} + 2|H|^{2} + 2\Delta logf ,
\end{equation}
where $A$ is the obstruction to integrability of the horizontal distribution and $H$ is the geodesic curvature of the fibers $S^{1}.$ Further $|H|^{2} = |\nabla logf|^{2},$ where log denotes the natural logarithm.

 Let $v(r)$ = area$D_{x_{o}}(r)$ in $(V, g_{V})$ and $\lambda (r) = v' (r)$ the length of $S_{x_{o}}(r),$ for some fixed base point $x_{o}\in V.$ The following general result uses only the submersion equation (4.2) and the assumption $s \geq $ 0. (The ${\cal R}_{s}^{2}$ equations are used only after this result).
\bbgin{proposition} \label{p 4.1. }
  Let (V, $g_{V})$ be as above, satisfying (4.1) with $s \geq $ 0. Then there exists a constant $c <  \infty $ such that
\begin{equation} \label{e4.2}
v(r) \leq  c\cdot  r^{2}, 
\end{equation}

\begin{equation} \label{e4.3}
\lambda (r) \leq  c\cdot  r, 
\end{equation}
and
\begin{equation} \label{e4.4}
\int_{V}|\nabla logf|^{2}dA_{V} \leq  c, \ \  \int_{V}|A|^{2}dA_{V} \leq  c. 
\end{equation}

\end{proposition}
{\bf Proof:}
 We first prove (4.2) by combining (4.1) with the Gauss-Bonnet theorem. Thus, pair (4.1) with a smooth cutoff function $\eta^{2}$ where $\eta  = \eta (d),$ and $d$ is the distance function on $V$ from $x_{o}\in V.$ This gives
\begin{equation} \label{e4.5}
\int_{V}(\eta^{2}s + \eta^{2}|A|^{2} + 2\eta^{2}|\nabla logf|^{2} + 2\eta^{2}\Delta logf) = \int_{V}\eta^{2}s_{V}. 
\end{equation}
Now
\begin{equation} \label{e4.6}
2\int_{V}\eta^{2}\Delta logf = - 4\int_{V}\eta<\nabla\eta , \nabla logf> \  \geq  - 2\int_{V}\eta^{2}|\nabla logf|^{2} - 2\int_{V}|\nabla\eta|^{2}, 
\end{equation}
so that,
\begin{equation} \label{e4.7}
\int_{V}\eta^{2}s_{V} \geq  - 2\int_{V}|\nabla\eta|^{2}, 
\end{equation}
since $s \geq $ 0. The Gauss-Bonnet theorem implies that in the sense of distributions on ${\Bbb R} ,$ 
\begin{equation} \label{e4.8}
s_{V}(r) = 2[\chi (r) -  \kappa (r)]' , 
\end{equation}
where 
\begin{equation} \label{e4.9}
\chi (r) = \chi (B(r)), \ s_{V}(r) = \int_{S_{x_{o}}(r)}s_{V}, \ \kappa (r) = \int_{S_{x_{o}}(r)}\kappa  ,
\end{equation}
and $\kappa $ is the geodesic curvature of $S_{x_{o}}(r),$ c.f. [GL, Theorem 8.11]. Since $\eta  = \eta (d),$ by expressing (4.7) in terms of integration over the levels of $d$, one obtains, after integrating by parts,
\begin{equation} \label{e4.10}
\int_{0}^{\infty}(\eta^{2})' [\chi (r) -  \kappa (r)] \leq  \int_{0}^{\infty}(\eta' )^{2}\lambda (r). 
\end{equation}
Now choose $\eta  = \eta (r) = -\frac{1}{R}r +$ 1, $\eta (r) =$ 0, for $r \geq  R$. In particular, $(\eta^{2})'  \leq $ 0. It is classical, c.f. again [GL, p. 391], that
\begin{equation} \label{e4.11}
\lambda' (r) \leq  \kappa (r) .
\end{equation}
Then (4.10) gives
$$2\int_{0}^{R}(-\frac{1}{R}r+1)\frac{1}{R}\lambda' (r) \leq  \frac{1}{R^{2}}\int_{0}^{R}\lambda (r)dr + 2\int (-\frac{1}{R}r+1)\frac{1}{R}\chi (r)dr. $$
But
$$2\int_{0}^{R}(-\frac{1}{R}r+1)\frac{1}{R}\lambda' (r) = \frac{2}{R}\int\lambda'  -  \frac{2}{R^{2}}\int r\lambda'  = \frac{2\lambda}{R} -  \frac{2}{R}\lambda  + \frac{2}{R^{2}}\int\lambda . $$
Thus
$$\frac{1}{R^{2}}\int_{0}^{R}\lambda (r)dr = 2\int (-\frac{1}{R}r+1)\frac{1}{R}\chi (r)dr \leq  C, $$
where the last inequality follows from the fact that $\chi (B(r)) \leq $ 1, since $B(r)$ is connected. This gives (4.2).

 Next, we prove (4.4). Returning to (4.5) and applying the Gauss-Bonnet theorem with $\eta  = \chi_{B(r)}$ the characteristic function of $B(r) = B_{x_{o}}(r)$ in (4.10) gives,
$$\int_{D(r)}(2|\nabla logf|^{2}+|A|^{2}) \leq  4\pi\chi (B(r)) -  2\lambda'-  2\int_{S(r)}<\nabla logf, \nu> , $$
where $\nu $ is the unit outward normal; the term $\lambda' $ is understood as a distribution on ${\Bbb R}^{+}.$ Using the H\"older inequality on the last term, we have
\begin{equation} \label{e4.12}
2\int_{D(r)}(|\nabla logf|^{2}+\frac{1}{2}|A|^{2}) \leq  4\pi\chi (B(r)) -  2\lambda' + 2(\int_{S(r)}|\nabla logf|^{2}+\frac{1}{2}|A|^{2})^{1/2}(v' )^{1/2}. 
\end{equation}
We now proceed more or less following a well-known argument in [CY,Thm.1]. Set 
$$F(r) = \int_{D(r)}(|\nabla logf|^{2}+\frac{1}{2}|A|^{2}). $$
Then (4.12) gives
\begin{equation} \label{e4.13}
F(r) \leq  2\pi  -  \lambda'  + (F' )^{1/2}\lambda^{1/2}. 
\end{equation}
Note that $F(r)$ is monotone increasing in $r$. If $F(r) \leq  2\pi ,$ for all $r$, then (4.4) is proved, so that we may assume that $F(r) >  2\pi ,$ for large $r$. Then (4.13) implies
$$\frac{1}{\lambda^{1/2}} + \frac{\lambda'}{(F- 2\pi )\lambda^{1/2}} \leq  \frac{(F' )^{1/2}}{F- 2\pi}. $$
We integrate this from $r$ to $s$, with $s >>  r$. A straightforward application of the H\"older inequality gives
$$(s- r)^{3/2}\cdot \bigl(\int_{r}^{s}\lambda\bigr)^{- 1/2} \leq  \int_{r}^{s}(\lambda^{- 1/2}). $$
Thus, using (4.2), one obtains
$$(s^{1/2} -  r^{1/2}) + \int_{r}^{s} \frac{\lambda'}{(F- 2\pi )\lambda^{1/2}} \leq  c\cdot  (F^{-1}(r) -  F^{-1}(s))^{1/2}(s-r)^{1/2}. $$
The integral on the left is non-negative, (integrate by parts to see this), and thus taking a limit as $s \rightarrow  \infty $ implies, for any $r$,
$$c \leq  F^{-1}(r), $$
which gives (4.4). 

 To prove (4.3), write (4.13) as $\lambda'  \leq  2\pi  + (F' )^{1/2}\cdot \lambda^{1/2}.$ Then one may integrate this from 0 to $r$ and apply the H\"older inequality and (4.2) to obtain (4.3).

{\endproof}

\noindent
\bbgin{remark} \label{r 4.2.}
  {\rm Consider the metric $\Roof{g}{\widetilde}_{V} = f^{2}\cdot  g_{V}$ on $V$. If $\Roof{K}{\widetilde}$ and $K$ denote the Gauss curvatures of $\Roof{g}{\widetilde}_{V}$ and $g_{V}$ respectively, then a standard formula, c.f. [B, 1.159], gives $f^{2}\Roof{K}{\widetilde} = K -  \Delta logf.$ Thus, from (4.1), we obtain
$$f^{2}\Roof{K}{\widetilde} = {\tfrac{1}{2}} s + {\tfrac{1}{2}} |A|^{2} + |\nabla logf|^{2} \geq  0. $$
If one knew that $\Roof{g}{\widetilde}_{V} $ were complete, then the results in Proposition 4.1 could be derived in a simpler way from the well-known geometry of complete surfaces of non-negative curvature.}

\end{remark}

 Next, we have the following analogue of Lemma 2.4, or more precisely (2.16). As in \S 2, let $t(x) = dist_{N}(x, x_{o}),$ for some fixed point $x_{o}\in N.$
\bbgin{proposition} \label{p 4.3.}
  Let (N, g, $\omega )$ be a complete ${\cal R}_{s}^{2}$ solution, with a free isometric $S^{1}$ action, and $\omega  \leq $ 0. Then there is a constant $c_{o} > $ 0 such that on (N, g)
\begin{equation} \label{e4.14}
\rho  \geq  c_{o}\cdot  t, 
\end{equation}
and
\begin{equation} \label{e4.15}
|r| \leq  c_{o}^{-1}t^{-2}. 
\end{equation}
Further, for any $x,y\in S(r),$
\begin{equation} \label{e4.16}
\frac{\omega (x)}{\omega (y)} \rightarrow  1, \ \ {\rm as} \ \ r \rightarrow  \infty , 
\end{equation}
and $\omega $ is a proper exhaustion function on N.

 If further $\omega $ is bounded below, then
\begin{equation} \label{e4.17}
lim_{t\rightarrow\infty}t^{2}|r| = 0, 
\end{equation}
and
\begin{equation} \label{e4.18}
osc_{N \setminus B(r)} \ \omega  \rightarrow  0, \ \ {\rm as} \ r \rightarrow  \infty , 
\end{equation}

\end{proposition}
{\bf Proof:}
 These estimates will be proved essentially at the same time, but we start with the proof of (4.14). This is proved by contradiction, and the proof is formally identical to the proof of (2.16). Thus, assuming (4.14) is false, let $\{x_{i}\}$ be any sequence in $(N, g)$ with $t(x_{i}) \rightarrow  \infty ,$ chosen to be $(\rho ,\frac{1}{2})$ buffered as following (2.16). Blow-down the metric $g$ based at $x_{i},$ by setting $g_{i} = \rho (x_{i})^{-2}\cdot  g.$ In the $g_{i}$ metric, the equations (0.4)-(0.5) take the form
\begin{equation} \label{e4.19}
\alpha_{i}\nabla{\cal R}^{2} + L^{*}\omega  = 0, 
\end{equation}

\begin{equation} \label{e4.20}
\Delta\omega  = -\frac{\alpha_{i}}{4}|r|^{2}, 
\end{equation}
where $\alpha_{i} = \alpha\rho_{i}^{-2}, \rho_{i} = \rho (x_{i}, g)$. All other metric quantities in (4.19)-(4.20) are w.r.t. $g_{i}.$ 

 For clarity, it is useful to separate the discussion into non-collapse and collapse cases.

{\bf Case (I).(Non-Collapse).}
 Suppose there is a constant $a_{o} > $ 0 such that
\begin{equation} \label{e4.21}
areaD(r) \geq  a_{o}r^{2}, 
\end{equation}
in $(V, g_{V}).$

 Now by Proposition 4.1,
\begin{equation} \label{e4.22}
\int_{V \setminus D(r)}|\nabla logf|^{2} \rightarrow  0, \ \ {\rm as} \ r \rightarrow  \infty . 
\end{equation}
The integral in (4.22) is scale-invariant, and also invariant under multiplicative renormalizations of $f$. Thus, we normalize $f$ at each $x_{i}$ so that $f(x_{i}) \sim $ 1. This is equivalent to passing to suitable covering or quotient spaces of $N$.

 By the discussion on convergence preceding Lemma 2.1, and by Theorem 3.6, a subsequence of the metrics $(B_{x_{i}}(\frac{3}{2}), g_{i}) \subset (N, g_i)$ converges smoothly to a limit ${\cal R}_{s}^{2}$ solution defined at least on $(B_{x}(\frac{5}{4}), g_{\infty})$, $x = lim x_i$, possibly with $\alpha = 0$, (i.e. a solution of the static vacuum Einstein equations (1.8)). In particular, $f$ restricted to $B_{x_{i}}(\frac{5}{4})$ is uniformly bounded, away from 0 and $\infty .$ It follows then from the bound on $f$ and (4.22) that
\begin{equation} \label{e4.23}
\int_{B_{x_{i}}(\frac{5}{4})}|\nabla logf|^{2} \rightarrow  0, \ \ {\rm as} \ \ r \rightarrow  \infty . 
\end{equation}
w.r.t. the $g_{i}$-metric on $N$.

 The smooth convergence as above, and (4.23) imply that on the limit, $f =$ const. The same reasoning on (4.4) implies that $A =0$ on the limit. Since $A = 0$, $\nabla f =$ 0 and $s =$ 0 on the limit, from (4.1) we see that $s_{V} =$ 0 on the limit. Hence the limit is a flat product metric on $C\times S^{1},$ where $C$ is a flat 2-manifold. The smooth convergence implies that $(B_{x_{i}}(\frac{5}{4}), g_{i})$ is almost flat, i.e. its curvature is almost 0. This is a contradiction to the buffered property of $x_{i},$ as in Lemma 2.4. This proves (4.14). The estimate (4.15) follows immediately from (4.14) and the smooth convergence as above.

 The proof of (4.17) is now again the same as in Lemma 2.4. Namely repeat the argument above on the metrics $g_{i} = t(x_{i})^{-2}\cdot  g,$ for any sequence $x_{i}$ with $t(x_{i}) \rightarrow  \infty .$ Note that in this (non-collapse) situation, the hypothesis that $\omega $ is bounded below is not necessary.

 Observe also that the smooth convergence and (4.22) imply that
\begin{equation} \label{e4.24}
|\nabla logf|(x) <<  1/t(x), 
\end{equation}
as $t(x) \rightarrow  \infty ,$ so that $f$ grows slower than any fixed positive power of $t$. It follows from this and from (4.3) that the annuli $A_{x_{o}}(r, 2r)$ have diameter satisfying diam$A_{x_{o}}(r,2r) \leq  c\cdot  r,$ for some fixed constant $c <  \infty .$

 To prove (4.16), we see from the above that for any sequence $x_{i}$ with $t(x_{i}) \rightarrow  \infty ,$ the blow-downs $g_{i}$ as above converge smoothly to a solution of the static vacuum Einstein equations (1.8), which is flat, (in a subsequence). Hence the limit potential $\bar{\omega} =$ lim $\bar{\omega}_{i}, \bar{\omega}_{i} = \omega (x)/|\omega (x_{i})|,$ is either constant or a non-constant affine function. The limit is a flat product $A(k^{-1},k)\times S^{1},$ where $A(k^{-1},k)$ is the limit of the blow-downs of $A(k^{-1}t(x_{i}),kt(x_{i}))\subset (V, g_{i})$ and $k > $ 0 is arbitrary. If $\bar{\omega}$ were a non-constant affine function, then $\bar{\omega}$ must assume both positive and negative values on $A(k^{-1}, k)$, for some choice of sufficiently large $k$, which contradicts the assumption that $\omega  \leq $ 0 everywhere. Thus, $\omega $ renormalized as above converges to a constant on all blow-downs, which gives (4.16).

 Since $N = {\Bbb R}^{2}\times S^{1}$ topologically, $N$ has only one end. By the minimum principle applied to the trace equation (3.2), $inf_{S(r)}\omega  \rightarrow  inf_{N}\omega ,$ as $r \rightarrow  \infty $. Together with (4.16), it follows that $\omega $ is a proper exhaustion function on $N$. Further, if $\omega $ is bounded below, then (4.18) follows immediately.

{\bf Case (II). (Collapse).}

 If (4.21) does not hold, so that
\begin{equation} \label{e4.25}
areaD(r_{i}) <<  r_{i}^{2}, 
\end{equation}
for some sequence $r_{i} \rightarrow  \infty ,$ then one needs to argue differently, since in this case, the estimate (4.22) may arise from collapse of the area, and not the behavior of $logf$.

 First we prove (4.14). By the same reasoning as above in Case (I), if (4.14) does not hold, then there is a $(\rho, \frac{1}{2})$ buffered sequence ${x_{i}}$, with $t(x_i) \rightarrow \infty$, which violates (4.14). Further, we may choose the base points exactly as in the proof of Lemmas 2.1 or 2.4, to satisfy (2.4). As in Case (I), normalize $f$ at each $x_i$ so that $f(x_i) \sim 1$. If the metrics $g_{i} = \rho(x_{i})^{-2} \cdot g$ are not collapsing at ${x_{i}}$, then the same argument as in Case (I) above gives a contradiction. Thus, assume the metrics $\{g_i\}$ are collapsing at $x_i$. Now as in (2.4), the curvature radius $\rho_i$ is uniformly bounded below by $\frac{1}{2}$ within arbitrary but fixed distances to the base point $x_{i}$. Hence, as discussed preceding Lemma 2.1, we may unwrap the collapse by passing to sufficiently large finite covers of arbitrarily large balls $B_{x_{i}}(R_i)$, $R_i \rightarrow \infty$. One thus obtains in the limit a complete non-flat ${\cal R}_{s}^{2}$ or static vacuum solution $(N' , g' )$, (corresponding to $\alpha = 0$), with an additional free isometric $S^{1}$ action, i.e. on $(N' , g' )$ one now has a free isometric $S^{1}\times S^{1}$ action. The second $S^{1}$ action arises from the collapse, and the unwrapping of the collapse in very large covers. This means that $N' $ is a torus bundle over ${\Bbb R} .$ Since $(N' , g' )$ is complete and scalar-flat, a result of Gromov-Lawson [GL, Thm. 8.4] states that any such metric is flat. This contradiction then implies (4.14) must hold. As before, smooth convergence then gives (4.15).

 The argument for (4.16)-(4.18) proceeds as follows. Consider any sequence $\{x_{i}\}$ in $N$ with $t(x_{i}) \rightarrow  \infty .$ The blow-down metrics $g_{i} = t(x_{i})^{-2}\cdot  g$ have a subsequence converging, after passing to suitably large finite covers as above, to a (now non-complete) maximal limit $(N' , g' )$ with, as above, a free isometric $S^{1}\times S^{1}$ action. The limit $(N' , g' )$ is necessarily a solution of the static vacuum equations (1.8) with potential $\bar \omega$ obtained by renormalizing the potential $\omega $ of $(N, g)$, i.e. $\bar{\omega} = lim_{i\rightarrow\infty} \frac{\omega (x)}{|\omega (x_{i})|}$ as before. Now it is standard, c.f. [An2, Ex.2.11], [EK, Thm.2-3.12], that the only, (even locally defined), solutions of these equations with such an $S^{1}\times S^{1}$ action are (submanifolds of) the {\it  Kasner metrics}, given explicitly as metrics on ${\Bbb R}^{+}\times S^{1}\times S^{1}$ by
\begin{equation} \label{e4.26}
dr^{2}+r^{2\alpha}d\theta_{1}^{2}+r^{2\beta}d\theta_{2}^{2}, 
\end{equation}
where $\alpha  = (a- 1)/(a- 1+a^{-1}), \beta  = (a^{-1}- 1)/(a- 1+a^{-1}),$ with potential $\bar{\omega} =  cr^{\gamma}, \gamma  = (a- 1+a^{-1})^{-1}$. The parameter $a$ may take any value in $[- 1,1].$ The values $a = 0$ and $a = 1$ give the flat metric, with $\bar{\omega} =$ const and $\bar{\omega} = cr$ respectively. 

  The limit $(N', g')$ is flat if $a = 0,1$; this occurs if and only if $|r|t^{2} \rightarrow $ 0 in a $g_i$-neighborhood of $x_{i}$. Similarly, the limit $(N', g')$ is non-flat when $a \in [-1, 0)\cup (0,1)$, which occurs when $|r|t^{2}$ does not converge to 0 everywhere in a $g_i$-neighborhood of $x_{i}.$ 

 In either case, the limits of the geodesic spheres in $(N, g)$ in the limit Kasner metric, are the tori $\{r =$ const\}. Since the oscillation of the limit potential $\bar{\omega}$ on such tori is 0, it is clear from the smooth convergence that (4.16) holds. As before, it is also clear that $\omega $ is a proper exhaustion function.

 Finally if $\omega $ is bounded below, then the limit $\bar \omega$ is uniformly bounded. Since $\bar \omega = cr^{\gamma}$, it follows that necessarily $a = 0$, so that all limits $(N' , g' )$ above are flat. The same argument as above in the non-collapse case then implies (4.17)-(4.18).

{\endproof}

\noindent
\bbgin{remark} \label{r 4.4.}
  {\rm The argument in Proposition 4.3 shows that if $\omega  \leq $ 0, then either the curvature decays faster than quadratically, i.e. (4.17) holds, or the ${\cal R}_{s}^{2}$ solution is asymptotic to the Kasner metric (4.26). It may be possible that there are complete $S^{1}$ invariant ${\cal R}_{s}^{2}$ solutions asymptotic to the Kasner metric, although this possibility remains unknown.}

\end{remark}

 We are now in position to complete the proof of Theorem 0.2.

\noindent
{\bf Proof of Theorem 0.2.}

 The assumption in (0.6) that $\omega $ is bounded below will be used only in one place, c.f. the paragraph following (4.37), so for the moment, we proceed without this assumption.

 It is convenient for notation to change sign, so we let
\begin{equation} \label{e4.27}
u = -\omega . 
\end{equation}
It is clear that Proposition 4.3 implies that $u > $ 0 outside a compact set in $N$.

 It is useful to consider the auxilliary metric
\begin{equation} \label{e4.28}
\Roof{g}{\widetilde} = u^{2}\cdot  g, 
\end{equation}
compare with [An2, \S 3]. In fact the remainder of the proof follows closely the proof of [An2, Thm.0.3], c.f. also [An2, Rmk.3.6] so we refer there for some details. We only consider $\Roof{g}{\widetilde}$ outside a compact set $K \subset  N$ on which $u > $ 0, so that $\Roof{g}{\widetilde}$ is a smooth Riemannian metric. By standard formulas for conformal change of the metric, c.f. [B, Ch.1J], the Ricci curvature $\Roof{r}{\widetilde}$ of $\Roof{g}{\widetilde}$ is given by
$$\Roof{r}{\widetilde} = r -  u^{-1}D^{2}u -  u^{-1}\Delta u\cdot  g + 2(dlog u)^{2} = $$

\begin{equation} \label{e4.29}
= - u^{-1}L^{*}u - 2u^{-1}\Delta u\cdot  g + 2(dlog u)^{2}. 
\end{equation}

$$\geq  - u^{-1}L^{*}u - 2u^{-1}\Delta u\cdot  g. $$
From the Euler-Lagrange equations (0.4)-(0.5) and (4.15), together with the regularity estimates from Theorem 3.6, we see that
$$|\nabla{\cal R}^{2}| \leq  c\cdot  t^{-4}, |\Delta u| \leq  c\cdot  t^{-4}. $$
Thus the Ricci curvature of $\Roof{g}{\widetilde}$ is almost non-negative outside a compact set, in the sense that
$$\Roof{r}{\widetilde} \geq  - c\cdot  t^{-4}\Roof{g}{\widetilde}. $$
Further, since the Ricci curvature controls the full curvature in dimension 3, one sees from (4.29) that the sectional curvature $\Roof{K}{\widetilde}$ of $\Roof{g}{\widetilde}$ satisfies
\begin{equation} \label{e4.30}
|\Roof{K}{\widetilde}| \leq  \frac{c}{u^{2}}|\nabla logu|^{2} + \frac{c}{t^4}, 
\end{equation}
where the norm and gradient on the right are w.r.t. the $g$ metric. 

 Let $\Roof{t}{\widetilde}(x) = dist_{\Roof{g}{\widetilde}}(x, x_{o}),$ (for $\Roof{t}{\widetilde}$ large), and $|\Roof{K}{\widetilde}|(\Roof{t}{\widetilde}) = sup_{S(\Roof{t}{\widetilde})}|\Roof{K}{\widetilde}|,$ taken w.r.t. $(N, \Roof{g}{\widetilde}).$ It follows from the change of variables formula that
\begin{equation} \label{e4.31}
\int_{1}^{\Roof{s}{\tilde}}\Roof{t}{\tilde}|\Roof{K}{\tilde}|(\Roof{t}{\tilde})d\Roof{t}{\tilde} \leq  c \bigl [\int_{1}^{s}t|\nabla logu|^{2}(t)dt + 1 \bigr ], 
\end{equation}
where as above, $|\nabla logu|^{2}(t) = sup_{S(t)}|\nabla logu|^{2}.$ In establishing (4.31), we use the fact that, 
\begin{equation} \label{e4.32}
d\Roof{t}{\tilde} = udt 
\end{equation}
together with the fact that
\begin{equation} \label{e4.33}
\Roof{t}{\tilde} \leq  c\cdot  u\cdot  t, 
\end{equation}
which follows from (4.32) by integration, using (4.16) together with the fact that $|\nabla logu| \leq  c/t,$ which follows from (4.15).

 Now we claim that
\begin{equation} \label{e4.34}
\int_{1}^{s}t|\nabla logu|^{2}(t)dt \leq  c\int_{1}^{s}area(S(t))^{-1}dt, 
\end{equation}
for some $c <  \infty .$ We refer to [An2,Lemma 3.5] for the details of this (quite standard) argument, and just sketch the ideas involved. First from the Bochner-Lichnerowicz formula and (4.15), one obtains 
$$\Delta|\nabla logu| \geq  - (c/t^{2})|\nabla logu|. $$
Hence, from the sub-mean value inequality, [GT, Thm.8.17], one has
\begin{equation} \label{e4.35}
sup_{S(r)}|\nabla logu|^{2} \leq  \frac{C}{volA(\frac{1}{2}r, 2r)}\int_{A(\frac{1}{2}r,2r)}|\nabla logu|^{2} \leq  \frac{C}{r\cdot  areaS(r)}\int_{B(r)}|\nabla logu|^{2}. 
\end{equation}
where the second inequality uses again the curvature bound (4.15). To estimate the $L^{2}$ norm of $|\nabla logu|$ on $N$, we observe that 
\begin{equation} \label{e4.36}
\int_{N}|r|^{2}dV = \int_{V}|r|^{2}fdA <  \infty . 
\end{equation}
The estimate (4.36) follows from the decay (4.15), from (4.2), and the fact that $sup_{S(r)}f \leq  r^{1+\varepsilon},$ for any fixed $\varepsilon  > $ 0, which follows from the proof of Proposition 4.3, (c.f. (4.24) and (4.26)). Now multiply the trace equation (0.5) by $u^{-1}$ and apply the divergence theorem on a suitable compact exhaustion $\Omega_{i}$ of $N$, for instance by sub-level sets of the proper exhaustion function $u$. Using (4.36), one thus obtains a uniform bound on the $L^{2}$ norm of $|\nabla logu|$ on $\Omega_{i},$ independent of $i$. Inserting this bound in (4.35) implies (4.34).

 Now if
\begin{equation} \label{e4.37}
\int_{1}^{\infty}area(S(t))^{-1}dt = \infty , 
\end{equation}
then a result of Varopoulos [V] implies that $(N, g)$ is parabolic, in the sense that any positive superharmonic function on $N$ is constant. 

 Since by (0.5) and (0.6), $\omega $ is a bounded superharmonic function on $(N, g)$, $\omega $ must be constant, so that $(N, g)$ is flat, by the trace equation (0.5). This proves Theorem 0.2 in this case. As indicated above, this is in fact the only place in the proof of Theorem 0.2 where the lower bound assumption $\omega  \geq  -\lambda  >  -\infty $ is used.

 Thus, suppose instead that
\begin{equation} \label{e4.38}
\int_{1}^{\infty}area(S(t))^{-1}dt <  \infty . 
\end{equation}
It follows then from (4.31),(4.34) and (4.38) that
\begin{equation} \label{e4.39}
\int_{1}^{\infty}\Roof{t}{\tilde}|\Roof{K}{\tilde}|(\Roof{t}{\tilde})d\Roof{t}{\tilde} <  \infty . 
\end{equation}

 Now it is a standard fact in comparison theory, c.f. [Ab], that the bound (4.39) implies that $(N, \Roof{g}{\widetilde})$ is almost flat at infinity, in the strong sense that geodesic rays starting at some base point $x_{o} \in N$ either stay a bounded distance apart, or grow apart linearly. More precisely, outside a sufficiently large compact set, $(N, \Roof{g}{\widetilde})$ is quasi-isometric to the complement of a compact set in a complete flat manifold. Observe that $(N, \Roof{g}{\widetilde})$ cannot be quasi-isometric to ${\Bbb R}^{3}$ outside a compact set, since that would imply that $V$ is non-collapsing at infinity. But then by combining (4.24), (4.34) and (4.38) it follows that $u\cdot  f,$ the length of the $S^{1}$ fiber in $(N, \Roof{g}{\widetilde})$ has sublinear growth; this is impossible when $(N, \widetilde g)$ is quasi-isometric to ${\Bbb R}^3$. Hence, outside a compact set, $(N, \widetilde g)$ is quasi-isometric to a flat product $C\times S^{1}$ where $C$ is the complement of a compact set in a complete flat 2-manifold.

 This means that there is a constant $C <  \infty $ such that the $\widetilde g$-length of the $S^{1}$ fiber satisfies
\begin{equation} \label{e4.40}
L_{\tilde g}(S^{1}) = u\cdot  f \leq  C. 
\end{equation}

 Now return to the trace equation (0.5) on $(N, g)$. Integrate this over $B(s) \subset (N, g)$ and apply the divergence theorem to obtain
\begin{equation} \label{e4.41}
\frac{\alpha}{4}\int_{B(s)}|r|^{2} = \int_{S(s)}<\nabla u, \nu> \  \leq\int_{S(s)}|\nabla u| = \int_{S_{V}(s)}|\nabla u|f, 
\end{equation}
where $S_{V}(s)$ is the geodesic $s$-sphere in $(V, g_{V}).$ Using (4.40), it follows that
$$\frac{\alpha}{4}\int_{B(s)}|r|^{2}  \leq  C\int_{S_{V}(s)}|\nabla logu|. $$
However, by (4.34) and (4.38), $|\nabla logu|(t) <<  t^{-1}.$ This and the length estimate (4.3) imply that
$$\int_{S_{V}(s)}|\nabla logu| \rightarrow  0, \ \ {\rm as} \ \  s \rightarrow  \infty , $$
which of course implies that $(N, g)$ is flat. This completes the proof of Theorem 0.2.

{\endproof}

\noindent
\bbgin{remark} \label{r 4.5.(i).}
  {\rm As noted above, the lower bound assumption on $\omega $ is required only in case $(N, g)$ is parabolic. Alternately, if $(V, g_{V})$ is non-collapsing at infinity, i.e. (4.21) holds, or if $u = - \omega$ has sufficiently small growth at infinity, i.e. $\int^{\infty}t|\nabla logu|^{2}dt <  \infty ,$ then the proof above shows that the lower bound on $\omega $ is again not necessary.

 On the other hand, as noted in Remark 4.4, there might exist complete ${\cal R}_{s}^{2}$ solutions asymptotic to the Kasner metric at infinity, so that $\omega  \sim  - r^{\gamma}, \gamma\in (0,1).$

{\bf (ii).}
 The proof of both Theorems 0.1 and 0.2 above only involve the asymptotic properties of the solution. It is clear that these proofs remains valid if it is assumed that $s \geq $ 0 outside a compact set in $(N, g)$ for Theorem 0.1, while $\omega  \leq $ 0 outside a compact set in $(N, g)$ for Theorem 0.2.}

\end{remark}

\section{Existence of Complete ${\cal R}_{s}^{2}$ Solutions.}
\setcounter{equation}{0}

{\bf \S 5.1.}
In this section, we show that the assumption that $(N, g, \omega )$ have an isometric free $S^{1}$ action in Theorem 0.2 is necessary, by constructing non-trivial ${\cal R}_{s}^{2}$ solutions with a large degree of symmetry.

 Let $g_{S}$ be the Schwarzschild metric on $[2m,\infty )\times S^{2},$ given by
\begin{equation} \label{e5.1}
g_{S} = (1-\frac{2m}{r})^{-1}dr^{2} + r^{2}ds^{2}_{S^{2}}. 
\end{equation}
The parameter $m > $ 0 is the mass of $g_{S}.$ Varying $m$ corresponds to changing the metric by a homothety. Clearly the metric is spherically symmetric, and so admits an isometric $SO(3)$ action, although the action of any $S^{1}\subset SO(3)$ on $[2m,\infty )\times S^{2}$ is not free.

 The boundary $\Sigma  = r^{-1}(2m)$ is a totally geodesic round 2-sphere, of radius $2m$, and hence $g_{S}$ may be isometrically doubled across $\Sigma $ to a complete smooth metric on $N = {\Bbb R}\times S^{2}.$

 The Schwarzschild metric is the most important solution of the static vacuum Einstein equations (1.8). The potential is given by the function
\begin{equation} \label{e5.2}
u = (1-\frac{2m}{r})^{1/2}, 
\end{equation}
The potential $u$ extends past $\Sigma $ as an {\it odd}  harmonic function under reflection in $\Sigma .$

 We show that there is a potential function $\omega $ on $N$ such that $(N, g_{S}, \omega )$ is a complete solution of the ${\cal R}_{s}^{2}$ equations (0.4)-(0.5) with non-zero $\alpha .$
\bbgin{proposition} \label{p 5.1.}
  The Schwarzschild metric (N, $g_{S})$ satisfies the ${\cal R}_{s}^{2}$ equations (0.4)-(0.5), where the potential is given by
$$\omega  = \tau  + c\cdot  u, $$
for any $c\in{\Bbb R} $ and $u$ is as in (5.2). The function $\tau $ is spherically symmetric and even w.r.t. reflection across $\Sigma .$ Explicitly,
\begin{equation} \label{e5.3}
\tau  = lim_{a\rightarrow 2m} \ \tau_{a} 
\end{equation}
where a $>  2m$ and
\begin{equation} \label{e5.4}
\tau_{a}(r) = \frac{\alpha}{8}(1-\frac{2m}{r})^{1/2}\bigl(\int_{a}^{r} \frac{1}{s^{5}(1-\frac{2m}{s})^{3/2}}ds -  \frac{1}{ma^{3}}\frac{1}{(1-\frac{2m}{a})^{1/2}}\bigr) , 
\end{equation}
for $r \geq  a$. In particular, $\tau (2m) = -\frac{\alpha}{8m}(2m)^{-3}, \tau' (2m) =$ 0, $\tau  < $ 0 everywhere, and $\tau $ is asymptotic to a constant $\tau_{o} < $ 0 at each end of $N$.

\end{proposition}
{\bf Proof:}
 By scaling, we may assume $m = \frac{1}{2}.$ One may rewrite the expression (0.2) for $\nabla{\cal R}^{2}$ via a standard Weitzenbock formula, c.f. [B, 4.71], to obtain
\begin{equation} \label{e5.5}
\nabla{\cal R}^{2} = \frac{1}{2}\delta dr + \frac{1}{2}D^{2}s - r\circ r -  R\circ r -  \frac{1}{2}\Delta s\cdot  g + \frac{1}{2}|r|^{2}\cdot  g. 
\end{equation}
The Schwarzschild metric $g_{S},$ or any spherically symmetric metric, is conformally flat, so that $d(r-\frac{s}{4}g) =$ 0. Since $g_{S}$ is scalar-flat, one thus has
\begin{equation} \label{e5.6}
\nabla{\cal R}^{2} = - r\circ r -  R\circ r + \frac{1}{2}|r|^{2}\cdot  g. 
\end{equation}
Let $t$ denote the distance to the event horizon $\Sigma $ and let $e_{i}, i =$ 1,2,3 be a local orthonormal frame on $N = S^{2}\times {\Bbb R} ,$ with $e_{3} = \nabla t,$ so that $e_{1}$ and $e_{2}$ are tangent to the spheres $t =$ const. Any such framing diagonalizes $r$ and $\nabla{\cal R}^{2}.$ The Ricci curvature of $g_{S}$ satisfies
$$r_{33} = - r^{-3}, r_{11} = r_{22} = \frac{1}{2}r^{-3}. $$
A straightforward computation of the curvature terms in (5.6) then gives
$$(\nabla{\cal R}^{2})_{33} = \frac{1}{4}r^{-6}, (\nabla{\cal R}^{2})_{11} = (\nabla{\cal R}^{2})_{22} = -\frac{1}{2}r^{-6}  $$

 We look for a solution of (0.4) with $\tau $ spherically symmetric, i.e. $\tau  = \tau (t).$ Then
$$D^{2}\tau  = \tau' D^{2}t + \tau'' dt\otimes dt,  \Delta\tau  = \tau' H + \tau'' , $$
where $H = \Delta t$ is the mean curvature of the spheres $t =$ const and the derivatives are w.r.t. $t$. One has $(D^{2}t)_{33} =$ 0 while $(D^{2}t)_{ii} = \frac{1}{2}H$ in tangential directions $i =$ 1,2. Thus, in the $\nabla t$ direction, (0.4) requires
\begin{equation} \label{e5.7}
\frac{\alpha}{4}r^{-6}+ \tau''  - (\tau' H+\tau'' ) -  \tau (- r^{-3}) = 0, 
\end{equation}
while in the tangential directions, (0.4) is equivalent to
\begin{equation} \label{e5.8}
-\frac{\alpha}{2}r^{-6}+ \frac{H}{2}\tau'  - (\tau' H+\tau'' ) -  \frac{\tau}{2}r^{-3} = 0. 
\end{equation}
The equations (5.7)-(5.8) simplify to the system
\begin{equation} \label{e5.9}
\tau' H -  \tau r^{-3} = \frac{\alpha}{4}r^{-6}, 
\end{equation}

\begin{equation} \label{e5.10}
\tau'' + \frac{H}{2}\tau'  + \frac{\tau}{2}r^{-3} = -\frac{\alpha}{2}r^{-6}. 
\end{equation}

 It is easily verified that (5.10) follows from (5.9) by differentiation, so only (5.9) needs to be satisfied. It is convenient to change the derivatives w.r.t. $t$ above to derivatives w.r.t. $r$, using the relation $\frac{dr}{dt} = (1-\frac{1}{r})^{1/2}.$ Since also $H = 2\frac{r'}{r} = \frac{2}{r}(1-\frac{1}{r})^{1/2},$ (5.9) becomes
\begin{equation} \label{e5.11}
\frac{2}{r}(1-\frac{1}{r})\dot {\tau} -  \tau r^{-3} = \frac{\alpha}{4}r^{-6}, 
\end{equation}
where $\dot {\tau} = d\tau /dr.$

 The linear ODE (5.11) in $\tau $ can easily be explicitly integrated and one may verify that $\tau $ in (5.3)-(5.4) is the unique solution on $([2m,\infty )\times S^{2}, g_{S})$ of (0.4)-(0.5) satisfying, (since $m = \frac{1}{2}),$
\begin{equation} \label{e5.12}
\tau (1) = -\frac{\alpha}{4} \ \ {\rm and} \ \ \frac{d}{dr}\tau|_{r=1} = 0. 
\end{equation}
It follows that the {\it even} reflection of $\tau $ across $\{r =$ 1\} $= \Sigma $ gives a smooth solution of (0.4) on the (doubled) Schwarzschild metric $(N, g_{S}),$ satisfying the stated properties. 

 Since Ker $L^{*} = \ <u> $ on the Schwarzschild metric, the potential $\omega $ has in general the form $\omega  = \tau +cu,$ for some $c\in{\Bbb R} .$

{\endproof}

 Thus the Schwarzschild metric, with potential function $\omega ,$ gives the simplest non-trivial solution to the ${\cal R}_{s}^{2}$ equations (0.4)-(0.5), just as it is the canonical solution of the static vacuum equations (1.8).

\noindent
\bbgin{remark} \label{r 5.2.}
  {\rm It is useful to consider some global aspects of solutions to the static vacuum equations in this context. It is proved in [An3,Appendix] that there are no non-flat complete solutions to the static vacuum equations (1.8) with potential $\omega  < $ 0 or $\omega  > $ 0 everywhere. Thus, this result is analogous to Theorem 0.1. While there do exist non-trivial complete smooth solutions with $\omega $ changing sign, for example the (isometrically doubled) Schwarzschild metric with $\omega  = \pm u$ in (5.2), such solutions are very special since they have a smooth event horizon $\Sigma  = \{\omega  =$ 0\}, c.f. [An2] for further discussion.

 Hence, among the three classes of equations considered here, namely the ${\cal R}^{2}$ equations (0.2), the static vacuum equations (1.8) and and the ${\cal R}_{s}^{2}$ equations (1.9), only the scalar curvature constrained ${\cal R}_{s}^{2}$ equations (1.9) for $0<\alpha<\infty $ admit non-trivial complete solutions with non-vanishing potential.

 In [An4], we will investigate in greater detail the structure of complete ${\cal R}_{s}^{2}$ solutions, 0 $<  \alpha  <  \infty ,$ as well as the structure of junction solutions, i.e. metrics which are solutions of the ${\cal R}^{2}$ equations in one region, and solutions of the ${\cal R}_{s}^{2}$ equations or static vacuum equations in another region, as mentioned in \S 1.}

\end{remark}

{\bf \S 5.2.}
 For the work in [An4] and that following it, it turns out to be very useful to have the results of this paper, and those of [AnI, \S 3], with the $L^{2}$ norm of the traceless Ricci curvature $z = r -  \frac{s}{3}g$ in place of $r$, i.e. with the functional 
\begin{equation} \label{e5.13}
{\cal Z}^{2} = \int|z|^{2}dV 
\end{equation}
in place of ${\cal R}^{2}.$ We show below that this can be done with relatively minor changes in the arguments.

 First, the Euler-Lagrange equations for ${\cal Z}^{2}$ are
\begin{equation} \label{e5.14}
\nabla{\cal Z}^{2} = D^{*}Dz + \frac{1}{3}D^{2}s -  2\stackrel{\circ}{R}\circ z + \frac{1}{2}(|z|^{2} -  \frac{1}{3}\Delta s)\cdot  g = 0, 
\end{equation}
\begin{equation} \label{e5.15}
tr\nabla{\cal Z}^{2} = -\frac{1}{6}\Delta s -  \frac{1}{2}|z|^{2} = 0. 
\end{equation}
Similarly, the ${\cal Z}_{s}^{2}$ equations for scalar-flat metrics, i.e. the analogues of (0.4)-(0.5), are
\begin{equation} \label{e5.16}
\alpha\nabla{\cal Z}^{2} + L^{*}(\omega ) = 0, 
\end{equation}
$$\Delta \omega  = -\frac{\alpha}{4}|z|^{2}. $$

 We begin by examining the functional $I_{\varepsilon}^{~-}$ given by
\begin{equation} \label{e5.17}
I_{\varepsilon}^{~-} = \varepsilon v^{1/3}\int|z|^{2} + (v^{1/2}\int (s^{-})^{2})^{1/2}, 
\end{equation}
i.e. the ${\cal Z}^{2}$ analogue of the functional $I_{\varepsilon}' $ in (1.16). Recall here that from \S 1 that the behavior of $I_{\varepsilon}' $ as $\varepsilon  \rightarrow $ 0 gives rise to the ${\cal R}^{2}$ and ${\cal R}_{s}^{2}$ equations. In the same way, the behavior of $I_{\varepsilon}^{~-}$ as $\varepsilon \rightarrow 0$ gives rise to the ${\cal Z}^2$ and ${\cal Z}_{s}^{2}$ equations.

 Again, as noted in \S 1, the existence and basic properties of critical points of the functional $I_{\varepsilon}$ in (1.1) were (essentially) treated in [An1, \S 8]. For this functional, the presence of $r$ or $z$ in the definition (1.1) makes no essential difference, due to the ${\cal S}^{2}$ term in (1.1). However, for the passage from ${\cal S}^{2}$ in (1.1) to ${\cal S}_{-}^{2}$ in (5.17), this is no longer the case, since we have no apriori control on $s^{+}$ = min$(s, 0)$. Thus, we first show that the results of [An1, \S 3 and \S 8] do in fact hold for metrics with bounds on $I_{\varepsilon}^{~-}$ and for minimizers of $I_{\varepsilon}^{~-}$.

 Let $r_{h}$ be the $L^{2,2}$ harmonic radius, as in [An1, Def. 3.2]. The main estimate we need is the following; compare with [An1, Rmk. 3.6].
\bbgin{proposition} \label{p 5.3.}
  Let $D$ be a domain in a complete Riemannian manifold $(N, g)$, such that
\begin{equation} \label{e5.18}
\varepsilon\int_{D}|z|^{2} + \int_{D}(s^{-})^{2} \leq  \Lambda , 
\end{equation}
where 0 $<  \varepsilon  <  \infty .$ Then there is a constant $r_{o} = r_{o}(\Lambda , \varepsilon ) > $ 0 such that
\begin{equation} \label{e5.19}
r_{h}(x) \geq  r_{o}\cdot \nu (x)\cdot \frac{dist(x, \partial D)}{diam D}. 
\end{equation}

\end{proposition}
{\bf Proof:}
 The proof is a modification of the proof of [An1, Thm. 3.5]. Thus, if (5.19) is false, then there exists a sequence of points $x_{i}$ in domains $(D_{i}, g_{i})$ such that
\begin{equation} \label{e5.20}
\frac{r_{h}(x_{i})}{dist(x_{i}, \partial D_{i}))} <<  \frac{\nu (x_{i})}{diam D_{i}} \leq  1, 
\end{equation}
where the last inequality follows since the volume radius $\nu$ is at most the diameter. Choose the points $x_{i}$ to realize the mimimal value of the ratio in (5.20). It follows as in the proof of [An1, Thm. 3.5] that a subsequence of the rescaled pointed manifolds $(D_{i}, g_{i}' , x_{i}), g_{i}'  = r_{h}(x_{i})^{-2}\cdot  g_{i},$ converges in the weak $L^{2,2}$ topology to a complete, non-compact limit manifold $(N, g' , x)$, with $L^{2,2}$ metric $g' ,$ and of infinite volume radius. From (5.20), one easily deduces that $r_{h}(y_{i}) >  \frac{1}{2},$ for all $y_{i}$ such that $dist_{g_{i}'}(x_{i}, y_{i}) \leq  R$, for an arbitrary $R <  \infty ,$ and for $i$ sufficiently large.

 Now the bound (5.18) and scaling properties of curvature imply that
\begin{equation} \label{e5.21}
z \rightarrow  0, \  \ s^{-} \rightarrow  0, 
\end{equation}
strongly in $L^{2},$ uniformly on compact subsets of $(N, g' , x)$. Hence the limit $(N, g' )$ is of constant curvature, and non-negative scalar curvature. In particular, the scalar curvature $s' $ of $g' $ is constant. If $s'  > $ 0, then $(N, g' )$ must be compact, which is impossible. Hence $s'  =$ 0, and so $(N, g' )$ is flat. Hence, $(N, g' )$ is isometric to the complete flat metric ${\Bbb R}^{3},$ since $(N, g' )$ has infinite volume radius.

 We claim that $s \rightarrow $ 0 strongly in $L^{2}$ also. Together with (5.21), this implies $r \rightarrow $ 0 strongly in $L^{2},$ and the proof proceeds as in [An1, Thm. 3.5]. To prove the claim, let $\phi$ be any function of compact support in $B_{y}(\frac{1}{2}),$ for an arbitrary $y\in N.$ With respect to the $g_{i}' $ metric, and by use of the Bianchi identity, we have
$$\int s\cdot \Delta \phi = -\int<\nabla s, \nabla \phi>  = \tfrac{1}{6} \int<\delta z, \nabla \phi>  = \tfrac{1}{6} \int< z, D^{2}\phi> , $$
where we have removed the prime and subscript $i$ from the notation. By $L^{2}$ elliptic regularity w.r.t. the metrics $g_{i}' ,$ (which are uniformly bounded locally in $L^{2,2}),$ it follows that
\begin{equation} \label{e5.22}
\int_{B_{y}(\frac{1}{4})}(s -  s_{o})^{2} \leq  C\cdot \int_{B_{y}(\frac{1}{2})}|z|^{2}, 
\end{equation}
where $s_{o}$ is the mean value of $s$ on $B_{y}(\frac{1}{2}),$ and where $C$ is a fixed constant, independent of $i$ and $y$. Thus, $s = s_{g_{i}}'$ converges strongly to its limit mean value in the limit. Since $s_{o} =$ 0 in the limit, it follows that $s \rightarrow $ 0 strongly in $L^{2},$ as required.

{\endproof}

  Let $\rho_{z}$ be the $L^{2}$ curvature radius w.r.t. $z$, i.e. again as in [An1, Def. 3.2], $\rho_{z}(x)$ is the largest radius of a geodesic ball about $x$ such that for $y\in B_{x}(\rho_{z}(x))$ and $D_{y}(s) = B_{x}(\rho_{z}(x))\cap B_{y}(s),$
\begin{equation} \label{e5.23}
\frac{s^{4}}{vol D_{y}(s)}\int_{D_{y}(s)}|z|^{2} \leq  c_{o}, 
\end{equation}
where $c_{o}$ is a fixed small constant. Of course $\rho_{r} \leq \rho_{z}$. Note that Proposition 5.3 implies in particular that for the $L^{2}$ curvature radius $\rho  = \rho_{r}$ as in \S 2,
\begin{equation} \label{e5.24}
\rho_{r}(x) \geq  \rho_{o}\cdot \rho_{z}(x), 
\end{equation}
where $\rho_{o}$ is a constant depending only a lower bound on $\nu $ and a bound on ${\cal S}  _{-}^{2}$ on $\rho_{z}(x).$ Note also that this local result is false without a local bound on ${\cal S}_{-}^{2}.$ For in this case, a metric of constant negative curvature, but with arbitrarily large scalar curvature will make $r_{h}$ arbitrarily small without any change to $\rho_{z}.$

 Proposition 5.3 shows that the analogue of [An1, Thm.3.5/Rmk.3.6] holds for the functional $I_{\varepsilon}^{~-}$ in place of ${\cal R}^{2},$ for any given $\varepsilon  > $ 0, as does [An1,Thm.3.7]. Given this local $L^{2,2}$ control, an examination of the proofs shows that the results [An1,Thm.3.9-Cor.3.11] also hold w.r.t. $I_{\varepsilon}^{~-},$ without any further changes, as does the main initial structure theorem, [An1, Thm.3.19].

 The local results [An1, Lem.3.12-Cor.3.17], dealing with collapse within the $L^{2}$ curvature radius will not hold for $\rho_{z}$ without a suitable local bound on the scalar curvature. However, as noted at the bottom of [An1,p.223], [An1, Lemmas 3.12,3.13] only require a lower bound on the $L^{2}$ norm the negative part of the Ricci curvature. Hence these results, as well as [An1, Cor.3.14-Cor.3.17] hold for metrics satisfying 
\begin{equation} \label{e5.25}
s \geq  -\lambda , 
\end{equation}
for some fixed constant $\lambda  >  -\infty .$ In particular, for metrics satisfying (5.25), we have
\begin{equation} \label{e5.26}
\rho_{r}(x) \geq  \rho_{o}\cdot \rho_{z}(x), 
\end{equation}
where $\rho_{o} = \rho_{o}(\lambda , c_{o})$ is independent of the volume radius $\nu .$

 We now use the results above to prove the following:
\bbgin{proposition} \label{p 5.4.}
  Theorems 0.1 and 0.2 remain valid for complete non-compact ${\cal Z}^{2}$ and ${\cal Z}_{s}^{2}$ solutions respectively.
\end{proposition}
{\bf Proof:}
 For Theorem 0.2, this is clear, since a scalar-flat ${\cal R}_{s}^{2}$ solution is the same as a scalar-flat ${\cal Z}_{s}^{2}$ solution.

 For Theorem 0.1, as noted in the beginning of \S 2, the form of the full Euler-Lagrange equations (0.2) or (5.14) makes no difference in the arguments, since elliptic regularity may be obtained from either one. Thus we need only consider the difference in the trace equations (0.3) and (5.15). Besides the insignificant difference in the constant factor, the only difference in these equations is that $r$ is replaced by $z$; the fact that the sign of the constants is the same is of course important.

 An examination of the proof shows that all arguments for ${\cal R}^{2}$ solutions remain valid for ${\cal Z}^{2}$ solutions, with $r$ replaced by $z$, except in the following two instances:

 (i). The passage from (2.7) to (2.8) in the proof of Proposition 2.2, which used the obvious estimate $|r|^{2} \geq  s^{2}/3.$ This estimate is no longer available for $|z|^{2}$ in place of $|r|^{2}.$

 (ii). In the proof of Lemma 2.9(ii), where $\Delta s(x_{i}) \rightarrow $ 0 implies $|r|(x_{i}) \rightarrow $ 0, and hence $s(x_{i}) \rightarrow $ 0, which again no longer follows trivially for $z$ in place of $r$.

 We first prove (ii) for ${\cal Z}^{2}$ solutions. By Lemma 2.1, we may assume that $(N, g)$ is a complete ${\cal Z}^{2}$ solution, with uniformly bounded curvature. Let $\{x_{i}\}$ be a minimizing sequence for $s \geq $ 0 on $(N, g)$. As before $\Delta s(x_{i}) \rightarrow $ 0, and so $|z|(x_{i}) \rightarrow $ 0, as $i \rightarrow  \infty .$ Since the curvature is uniformly bounded, it follows from (the proof of) the maximum principle, c.f. [GT, Thm.3.5], that $|z|^{2} \rightarrow $ 0 in balls $B_{x_{i}}(c),$ for any given $c <  \infty .$ Hence the metric $g$ in $B_{x_{i}}(c)$ approximates a constant curvature metric. Since $(N, g, x_{i})$ is complete and non-compact, this forces $s \rightarrow $ 0 in $B_{x_{i}}(c),$ which proves (ii).

 Regarding (i), it turns out, somewhat surprisingly, that the proof of Proposition 2.2 is not so simple to rectify. First, observe that Proposition 2.2 is used only in the following places in the proof of Theorem 0.1.

 (a). The end of Lemma 2.4.

 (b). Lemmas 2.8 and 2.10.

 Regarding (a), the estimate (2.19) still holds. In this case, as noted following the proof of Proposition 2.2, there is a smoothing $\tilde t$ of the distance function $t$ such that $|\Delta \tilde t| \leq  c/ \tilde t^{2}$ and in fact 
\begin{equation} \label{e5.27}
|D^{2}\tilde t| \leq  c/\tilde t^{2}. 
\end{equation}
The proof may then be completed as follows. As in the beginning of the proof of Proposition 2.2, we obtain from the trace equation (5.15),
\begin{equation} \label{e5.28}
\int\eta^{4}|z|^{2} \leq  \frac{1}{3}\int<\nabla s, \nabla \eta^{4}> = -\frac{1}{18}\int< z, D^{2}\eta^{4}> , 
\end{equation}
where the last equality follows from the Bianchi identity $\delta z = -\frac{1}{6}ds,$ and $\eta  = \eta (\tilde t)$ is of compact support. Expand $D^{2}\eta^{4}$ as before and apply the Cauchy-Schwarz inequality to (5.28). Using (5.27), together with the argument following (2.9), the proof of Proposition 2.2 follows easily in this case.

 (b).  For both of these results, it is assumed that (2.47) holds, i.e. $s(x) \geq  d_{o}\cdot \rho (x)^{-2},$ and so $s(x) \geq  d\cdot  t(x)^{-2}.$ In this case, the estimate (5.26), which holds since $(N, g)$ has non-negative scalar curvature, together with a standard covering argument implies that
\begin{equation} \label{e5.29}
\int_{B(R)}s^{2}dV \leq  c_{1}\cdot \int_{B(2R)}|z|^{2}, 
\end{equation}
for all $R$ large, where $c_{1}$ is a constant independent of $R$. As before in the proof, there exists a sequence $r_{i} \rightarrow \infty$ as $i \rightarrow \infty$ such that, for all $R\in [r_{i}, 10r_{i}],$ 
\begin{equation} \label{e5.30}
\int_{B(R)}s^{2}dV \leq  c_{2}\cdot \int_{B(R/2)}s^{2}, 
\end{equation}
with $c_{2}$ independent of $R$. Given the estimates (5.29)-(5.30), the proof of Proposition 2.2 then proceeds exactly as before.

 The remainder of the proof of Theorem 0.1 then holds for ${\cal Z}^{2}$ solutions; the only further change is to replace $r$ by $z$.

{\endproof}

\bbgin{remark} \label{r 5.5.}
  {\rm We take this opportunity to correct an error in [AnI]. Namely in [An1, Thms. 0.1/3.19], and also in [An1, Thms. 0.3/5.9], it is asserted that the maximal open set $\Omega $ is embedded in $M$. This assertion may be incorrect, and in any case its validity remains unknown. The error is in the statement that the diffeomorphisms $f_{i_{k}}$ constructed near the top of [AnI, p.229] can be chosen to be nested.

 The proof of these results does show that $\Omega $ is weakly embedded in $M$,$$\Omega  \subset\subset  M, $$
in the sense that any compact domain with smooth boundary in $\Omega $ embeds as such a domain in $M$. Similarly, there exist open sets $V \subset  M$, which contain a neighborhood of infinity of $\Omega ,$ such that $\{g_{i}\}$ partially collapses $V$ along F-structures. In particular, $V$ itself, as well as a neighborhood of infinity in $\Omega $ are graph manifolds. Thus, the basic structure of these results remains the same, provided one replaces the claim that $\Omega \subset M$ by the statement that $\Omega \subset\subset M$. 

 The remaining parts of these results hold without further changes. The same remarks hold with regard to the results of \S 8. My thanks to thank Yair Minsky for pointing out this error.}

\end{remark}

\input{ref.tex}
\input{tail.tex}

\end{document}

%% file: ref.tex
\bibliographystyle{plain}

%% file: tail.tex
\begin{center}
September, 1998/October, 1999
\end{center}
\medskip
\address{Department of Mathematics\\
S.U.N.Y. at Stony Brook\\
Stony Brook, N.Y. 11794-3651}\\
\email{anderson@@math.sunysb.edu}